\documentclass[preprint,12pt]{elsarticle}

\usepackage{amsmath}
\usepackage{amsthm}
\usepackage{amssymb}
\usepackage{amsfonts}
\usepackage{graphicx}
\usepackage{verbatim}
\usepackage{mathrsfs}
\usepackage{subfigure}
\usepackage{fancyhdr}
\usepackage{latexsym}
\usepackage{graphicx}
\usepackage{dsfont}

\theoremstyle{plain}
\newtheorem{theorem}{Theorem}[section]
\newtheorem{proposition}[theorem]{Proposition}
\newtheorem{lemma}[theorem]{Lemma}

\newtheorem{remark}[theorem]{Remark}

\theoremstyle{definition}

\newcommand{\w}{\omega}
\newcommand{\al}{\alpha}
\newcommand{\R}{\mathbb{R}}
\newcommand{\E}{\mathcal{E}}

\newcommand{\bydef}{:=} 

\newcommand{\vect}[1]{{\boldsymbol{#1}}}
\newcommand{\bfA}{\vect{A}}
\newcommand{\bfB}{\vect{B}}
\newcommand{\bfC}{\vect{C}}

\date{}

\begin{document}

\begin{frontmatter}

\title{ The monotonicity of the apsidal angle in power-law potential systems}
\author{Roberto Castelli\corref{Roberto Castelli}}

\address{Department of Mathematics,
Vrije Universiteit Amsterdam,
De Boelelaan 1081,
1081 HV Amsterdam,
The Netherlands; {\tt r.castelli@vu.nl}}

\begin{abstract}
In a central force system the apsidal angle is the angle at the centre of force between two consecutive apsis and measures the precession rate of the line of apsis. The apsidal angle has applications in different fields and 
the Newton's apsidal precession theorem  has been extensively studied by astronomers, physicist and mathematicians. The perihelion precession of Mercury, the dynamics of galaxies, the vortex dynamics,  the JWKB quantisation  condition are some examples where the apsidal angle is of interest.

In case of  eccentric orbits and forces far from inverse square, numerical investigations provide evidence of the monotonicity of the apsidal angle with respect to the orbit parameters, such as the orbit eccentricity. However, no proof of this statement is available.

In this paper central force systems with  $f(r)\sim\mu r^{-(\al+1)}$  are considered. We  prove that for any $-2<\alpha<1$ the apsidal angle is a monotonic function of the orbital eccentricity, or equivalently of the angular momentum. As a corollary, the conjecture stating the absence of isolated cases of zero precession is proved.

\end{abstract}
\end{frontmatter}
\section{Introduction}

In a centripetal force system, when the force depends only on the distance from the centre,  a particle oscillates between the points of highest and lower distance, called apsis,  while rotating around the centre of the force. The apsidal angle $\Delta\theta$ refers to the angle at the force spanned by the particle between two consecutive apsis. 

Since Newton, the behaviour of the apsidal angle has been  extensively studied, in particular for its implications to celestial mechanics. Newton's precession theorem, stated in the Book I of {\it Principia} asserts a relation between the magnitude of the force and the apsidal angle:  for close to circular orbits, a force proportional to $r^{n-3}$ leads to $\Delta\theta= \pi/\sqrt{n}$. By means of this formula, the experimental  measure of the apsidal precession results in the exponent in the force law of the underlying system. For instance, the quiescence of Moon's apsis and of  the planetary aphelia,  that is the non rotation of the planetary apsidal line, allowed Newton to conclude that the variation of the  gravitational attraction of the Sun and of the Earth is inverse square.

Although Newton derived his formula for almost circular orbits, it seems that he trusted the formula also for eccentric orbits, when he applied it to Mercury's and Halley comet's motion.   
The analysis of the Mercury's orbit  gave rise to a long standing problem and  fascinating  debate among astronomers. During the nineteenth century astronomers realised that the precession rate of the Mercury's apses does not agree with  Newton's  law. In 1859 Le Verrier  \cite{LeVerrier1859} and later in 1882 Newcomb \cite{Newcomb} respectively found that 38 and 43 arc seconds per century of the precession of Mercury's apses could not be accounted for  on the basis of the inverse square law. In order to explain the anomalous precession, Le Verrier suggested
the presence of Vulvano,  a inter-mercurian planet, never found. 
 Other theories, all rejected,   supposed the existence of  a satellite orbiting around Mercury or that the gap in precession is completely due to the oblateness of the Sun.
  In 1894 Hall \cite{Hall1894} postulated a different power law for the gravitation attraction. Starting from the Bertrand's formula \cite{Bertrand} $\Delta\theta=\pi/\sqrt{N+3}$, where $N=n-3$, Hall inferred that the right exponent that would account for the Mercury anomaly is $n=1-\delta$ with $\delta\sim 1.6\cdot 10^{-7}$.
 Like the others, Hall's hypothesis was not completely convincing and the mystery of  Mercury's precession was finally solved by Einstein and his theory of General Relativity. Nowadays, the perihelion precession of Mercury is considered one of the empirical proof of the theory of general relativity.
 Note that Hall used  Betrand's formula for sizeable eccentricity and non inverse-square potential forces whereas Bertrand,  
 like  Newton, inferred his formula under the hypothesis of small eccentricity.
 
 When $n-2\leq 0$ the potential function $V(r)$ giving the force $f=-\nabla V$ has a singularity at $r=0$ ( $V(r)\sim r^{n-2}$ when $n\neq 2$ and $V(r)\sim-\log(r)$ for $n=2$). In these cases, the limiting value of the apsidal angle for orbits approaching the collision determines whether the flow can be extended  beyond the collision \cite{Mcgehee, gbegfugfg, Stoica}. A regularisation of the singular flow is possible only if the apsidal angle tends to a multiple of $\pi/2$. For instance, the Levi Civita regularisation \cite{LeviCivita} of the Kepler problem is possible because the apsidal angle tends to $\pi$ as the angular momentum goes to zero. Similarly in the logarithmic potential problem, the regularisation obtained via the transmission solution is dictated by the fact that the apsidal angle tend to $\pi/2$ \cite{MR2836348, MR3057157}. Also, when dealing the $N$-body and the $N$-center problem with variational techniques, the knowledge  of the apsidal angle for close to collision solutions, may help in proving collision avoidance, see for instance \cite{MR2299759}.

 Besides   celestial mechanics related problems,  power law centripetal force systems arise as models in a wide range of fields. The logarithmic potential, for instance,  has applications to galactic dynamics \cite{Richstone1982, Schwa1979}, to vortex filament dynamics \cite{Newton, MR2029363} as well as in particle physics, potential theory and astrophysics, see \cite{MNR:MNR22071} and reference therein.
 The behaviour of the apsidal angle for homogenous potentials has implication in astrophysical, see for instance \cite{Valsecchi} and quantum mechanics contexts, as in describing  how the radial quantum number and the orbital angular momentum contribute in the JWKB quantisation condition \cite{Feldman1979}. 
 
 The determination of the apsidal angle is as well  an interesting mathematical problem of its own. As we will show, it involves solving an integral with singularities at both endpoints. Also, a perturbative approach to  the orbital  equation  reveals connections to   the Mathieu-Hill differential equation \cite{valluri99}. 
 
  Now we know that Newton was right  when he trusted his formula in the case of inverse law potential, but not for potentials different than $r^{-1}$. In  the latter case the apsidal angle generally depends on the  orbital parameter as, for instance,  the eccentricity.  Indeed, a celebrated theorem of Bertrand \cite{Bertrand} states that there are only two central fields for which all bounded orbits are closed, the harmonic oscillator and the inverse square. In other words, only when $n=1$ and $n=4$ the apsidal angle does not depend on the initial condition and it holds $\Delta\theta=\pi/\sqrt{n}$.
  
A natural question is how the apsidal angle behaves for force laws different from inverse square and linear as the orbital parameters change. In case of power law forces, that is the case we are interested in, the same question has been raised  by different authors and several  results are present in the literature.

Let us formulate  the one-centre problem as the problem of finding  solutions for
\begin{equation}\label{eq:system}
\ddot u=-\nabla V_\al(|u|),\quad  u\in \R^2,\quad {\rm where}\quad 
\left\{\begin{array}{ll}
{\displaystyle V_{\al}(x)=\frac{1}{\al}\frac{1}{x^{\al}},\quad }& \al\neq 0\\
\\
{\displaystyle V_{0}(x)=-\log(x) }&\al=0.
\end{array}\right.
\end{equation}
According to this notation, $\alpha=1$ corresponds to the inverse square force law, $\alpha=-2$ is the harmonic oscillator and $\alpha=0$ refers to the logarithmic potential problem.

The solutions of system \eqref{eq:system} preserve the mechanical energy $E=\tfrac{1}{2}|\dot u|^{2}-V$ and the angular momentum $\vec\ell=u\wedge \dot u$, hence the trajectories lie on the plane perpendicular to $\vec \ell$. 
For an admissible choice of $E$ and $\ell=|\vec \ell|$  for which the orbit $u(t)$  is bounded,  the apsidal angle   is given by
\begin{equation}\label{eq:intro}
\Delta_{\al}\theta=\int_{r_{-}}^{r_{+}}\frac{\ell}{r^{2}\sqrt{2\big(E+V_\al(r)\big)-\frac{\ell^{2}}{r^{2}}}}\, dr
\end{equation}
where $r_{\pm}$ are the distances of the apsis from the centre of force and correspond to the values   where the radicand  vanishes.  $\alpha$ being  fixed, the parameters are $E$ and $\ell$. However both $E$ and $\ell$ can be expressed as function of $r_\pm$, hence  the apsidal angle only depends on the value of the apsidal distances and, more precisely, on the ratio $r_+/r_-$. It is common to introduce the orbital eccentricity $e=(r_+-r_-)/(r_++r_-)$ as a unique parameter and to study the behaviour of $\Delta_\al\theta$ as function of $e$. The parameter $e$ coincides with the standard  eccentricity when the orbit is elliptical and it is well defined also for non elliptical orbits. Thus, it is a valuable  parameter in the families of bounded  solutions of the one centre problem for any $\al$.
Note that $e=0$ for circular orbit and $e=1$ for degenerate orbits, that is when the solution lies on a straight line and falls in the collision. The latter situation occurs if and only if  the angular momentum vanishes, while, for a choice of $E$,  the maximal allowed value of the angular momentum leads to the circular solution. Indeed $e$ is monotonically decreasing with respect to $\ell$.

The analysis of the apsidal angle can be restricted to $\al<2$. Indeed the case $\al=2$ corresponds to  the force $\sim \kappa r^{-3}$ and the solution is given by the Cotes spiral \cite{Danby} in which $u(t)$ spirals towards the centre of force from any initial condition. When $\al=2$ and, similarly, when $\al>2$ no bounded orbit exists. 
 
 Since the integrand is singular at both the end-points and the apsidal values  themselves are only implicitly defined as functions of $E$ and $\ell$, both the analytical and numerical  treatment of the integral \eqref{eq:intro} is a challenging problem. 
 In general, \eqref{eq:intro} can not be solved in closed form. Exceptional cases are $\al=1$ and $\al=-2$ for which 
  the dependence on $e$ disappears and the integral can be easily calculated providing $\Delta_1\theta=\pi$ and $\Delta_{-2}\theta=\pi/2$. When $\al=-1,-\frac{2}{3},\frac{1}{2},\frac{2}{3}$ the apsidal angle can be expressed in terms of the elliptic functions, see \cite{whitta}. For others values of $\al$ one may trust  numerical computation. However, as pointed out in \cite{MNR:MNR8819, 1997MNRAS},  a numerical integration  of $\Delta_\al\theta$ may be very time consuming and could lead to unreliable results, in particular when $\al\sim 1$.
 
Power series expansions of the apsidal angle have been proposed for instance in \cite{1997MNRAS} and \cite{MR2497199}. In the first one,  by means of the Mellin transform,  the apsidal angle is expanded in powers of $\ell$ for any value $\al\neq 0$. The power series then obtained has the property of being asymptotic when $\ell\to 0$ but generically not converging for fixed $\ell $. In that paper the authors also showed how the knowledge of the apsidal precession in the spherical symmetric potential can lead to the description of the orbits in case of non axisymmetric potentials.   In the second work a generalisation of the B$\ddot{\rm u}$rmann-Lagrange series for multivalued inverse function has been applied and the apsidal angle for any central force field   has been written as a series of powers of $E-U_0$, $U_0$ being  the minimum of the effective potential.  

 Analytic  approximation  have been employed to study the radial equation and the  apsidal angle for small eccentricity. In the logarithmic potential case Touma and Tremaine \cite{1997MNRAS} used the epyciclic approximation and found the value $\Delta_0\theta=\pi/\sqrt{2}$ in the limit of $e\to 0$, Struck in \cite{struck} and Valluri et al. in \cite{MNR:MNR22071} used $p$-ellipses curves and the Lambert W function to approximate the radial solution up to order $e^2$ and the apsidal angle to high accuracy.

By means of perturbative methods in \cite{MR1433938} the apsidal angle has been analytically studied for power law forces slightly differing from inverse square
 showing, for instance,  that $\Delta_{1+\delta}\theta\sim\pi(1+\frac{\delta}{2}(1+\frac{e^2}{4}+\frac{e^4}{8}))$ valid for $|\delta|<<1$. 
  In particular it has been inferred that, for $\delta>0$,  the apsidal angle is increasing when $e$ increases and the opposite when $\delta<0$.
  In the same papers, by means of numerical integration, the apsidal angle has been computed for $\al=\frac{1}{2},\frac{3}{4},\frac{5}{4},\frac{3}{2}$. In \cite{MNR:MNR8819} the authors extend the results of the previous work and analyse the  derivative of the apsidal with respect to the power law exponent $\al$ and the angular momentum $\ell$. Analytical expressions are produced in case  $\alpha$ close to $1$ while the cases $0<\al<2$ are  studied numerically. The logarithmic potential case is also treated.

The various numerical experiments produced  in the quoted works  seem to confirm the perturbative results obtained for small $\delta$: the apsidal angle is  decreasing in $e$ for $0\leq \al<1$. Also, it has been conjectured in \cite{MNR:MNR8819} that   zero precession never occurs for a  power law potential different from inverse square. To the author's knowledge, no proof of the monotonicity of the apsidal angle and of the absence of isolated cases of zero precession is available.

\begin{figure}[htbp]
\begin{center}
\includegraphics[width=\textwidth]{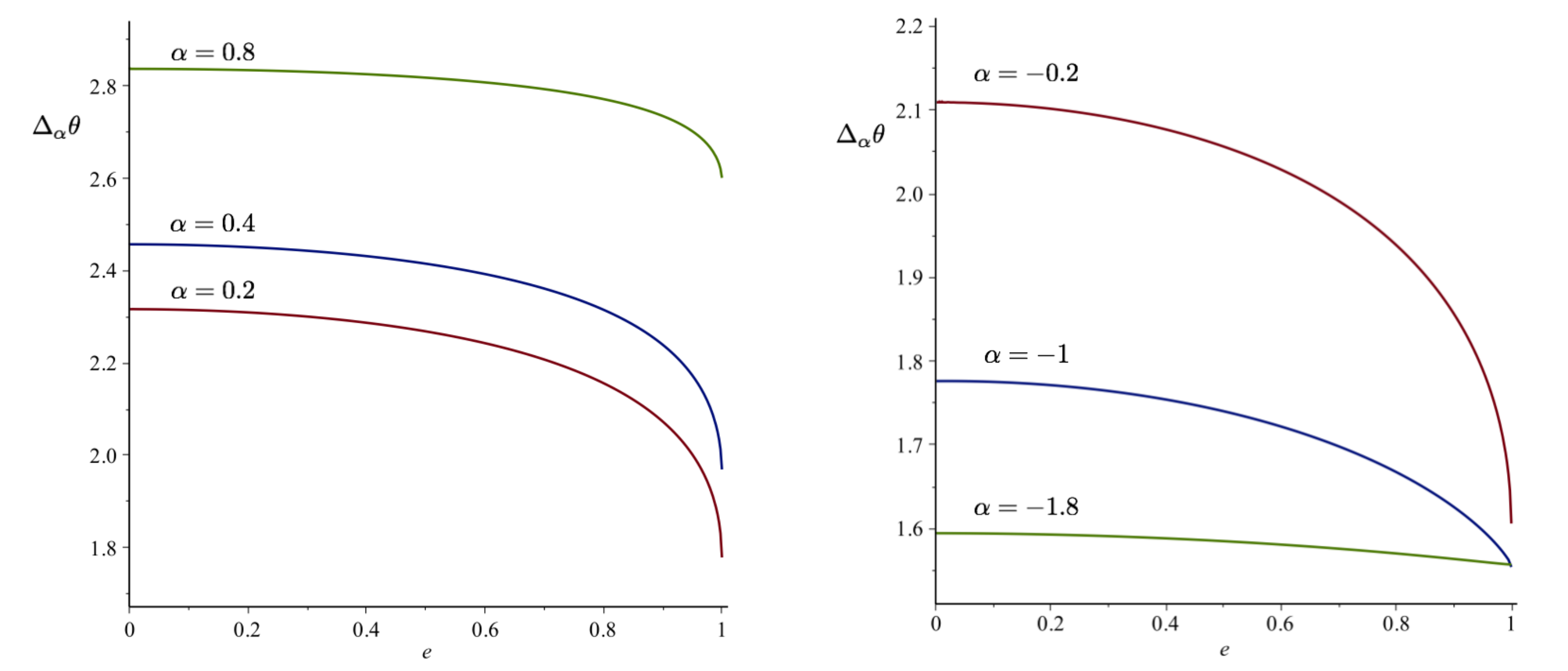}
\caption{Numerical evaluation of the apsidal angle as function of the eccentricity $e$ for different values of the power-law exponent $\alpha$.}
\label{fig:numerical}
\end{center}
\end{figure}

The aim of this paper is to  prove that for any $-2< \al<1$ the apsidal angle is a monotonic function of the orbital eccentricity  and, equivalently, of the value of the angular momentum $\ell$. In particular, as it appears from the numerical simulation reported in Figure  \ref{fig:numerical}, we prove that $\Delta_\al\theta$ is decreasing as a function of the eccentricity for any $\al\in(-2,1)$.  The limit values being already well known, the monotonicity allows to have a clear picture of the behaviour of the apsidal angle.  As a corollary, it also follows that the above conjecture is true. 

The monotonicity of the apsidal angle in the logarithmic potential case has been already proved in \cite{MR3159801}. We now extend that result for a homogenous potential $V_\al$ with $\al\neq 0$ and we prove the following.
\begin{theorem}\label{theo:main}
For any $\al\in(-2.1)$, the apsidal angle $\Delta_\al\theta$ is a monotonic function of the orbital eccentricity and 
$$\dfrac{\partial \Delta_\al\theta}{\partial e}<0.
$$
\end{theorem}

The proof of the theorem follows the same scheme adopted in \cite{MR3159801}. To begin with,  the integral \eqref{eq:intro} is recast   as a fixed endpoints integral and the contributions that depend on the orbital parameters are well separated. As already said, $\Delta_\al\theta$ depends on the ratio $r_+/r_-$ solely hence it can be written as a function of a unique parameter. Rather than using the orbital eccentricity $e$, we introduce a different parameter, denoted by $q$, that plays the same role as the eccentricity. 
The second step consists in the introduction of a function $\mathcal I_\al$ the positivity of which implies the monotonicity of the apsidal angle with respect to $q$. Such a function is then expanded  as a power series in the variable $q$. The third fundamental step consists in the analysis of the coefficients of the power series and in proving that the series is positive for any $\al,q$ of interest.

\section{The orbital equation and the apsidal angle}

Writing the trajectory $u(t)$ in polar coordinates, $u(t)=r(t)e^{i\theta(t)}$, the one centre problem \eqref{eq:system} is equivalent to solve the system
\begin{equation*}%\label{eq:systpol}
\ddot r -\frac{\ell^2}{r^3}+\frac{1}{r^{\al+1}}=0,\qquad \frac{d}{dt}(r^2\dot\theta)=0
\end{equation*}
where $\ell=r^2\dot\theta$ is the magnitude of the angular momentum $\vec\ell=u\wedge \dot u$.
The rotational symmetry of the  system implies the conservation of the angular moment and the second of the above equations simply formulates the conservation of  $\ell$.
Beside $\vec\ell$, the system \eqref{eq:system} conserves the total energy $E=\frac{1}{2}\dot r^{2}+\frac{1}{2}\frac{\ell^{2}}{r^{2}}-V_{\al}(r)$, from where we infer that the motion is allowed only for those values of the radial coordinate $r$ such that $E>\frac{1}{2}\frac{\ell^{2}}{r^{2}}-V_{\al}(r)$. The analysis of the potential function $V_\al$ leads to the following considerations: bounded orbits exist for negative values of $E$ if $\al\in (0,1)$, for all values of $E$ when $\al=0$ and for positive $E$ in case $\al\in(-2,0)$. Moreover the possible values for the  angular momentum $\ell$ are $\ell \in [0,\ell_{max}]$ where $\ell_{max}=\left( \frac{2\al}{\al-2}E\right)^{(\al-2)/2\al}$ if $\al\neq 0$ and $\ell_{max}=e^{E-1/2}$ for $\al=0$. For a selection of $E,\ell$ giving bounded solutions, the apsidal values $r_\pm$ are defined as zeros of the quantity $\frac{1}{2}\frac{\ell^{2}}{r^{2}}-V_{\al}(r)-E$.

In the logarithmic potential case, the apsidal equation assumes the form $\frac{\ell^2}{r^2}+\log(r^2)-2E=0$ the solutions of which can be expressed in terms of the two branches of the Lambert W function \cite{LambertW} as 
$$
\frac{1}{r}
=\sqrt{\frac{-W_{j}(-\ell^{2}e^{-2E})}{\ell}}.
$$
The relation $\ell=r^{2}\dot \theta$ implies that the rate of rotation between two consecutive apsis is given by $\Delta_\al\theta=\int_{r_{-}}^{r_{+}}\frac{\ell}{r^{2}\dot r}\, dr$, from where formula \eqref{eq:intro} comes. The apsidal angle is the main object of study of this paper and a  quick look at the integral \eqref{eq:intro} immediately reveals the main technical hurdles. Namely, the integral is singular at the endpoints and the endpoints are implicitly defined as solutions of the apsidal equation. We begin by reworking the apsidal angle into an integral with fixed endpoints.
To this goal it is convenient to reformulate the orbit equation in Clairaut's coordinates, where $\theta$ is the independent variable and $z=1/r$. In this setting the orbital differential equation is 
$$
\frac{d^{2}}{d\theta^{2}}z+z=f\left(\frac{1}{z}\right)\frac{1}{z^{2}\ell^{2}}
$$
 where $f(r)$ is the force.
We now  define $P(z)=\frac{1}{z^2}f(\frac{1}{z})$ and $w(z)=\int 2P(z)\, dz$. Then, according to  Bertrand \cite{Bertrand}, the apsidal angle can be written in the form
$$
\Delta_\alpha\theta=\int_a^b\frac{\sqrt{w(b)-w(a)}}
{\sqrt{a^2w(b)-b^2w(a)+(b^2-a^2)w(z)-z^2(w(b)-w(a))}}dz
$$
where $a$, $b$ are the apsidal values for $z$. Easy algebraic manipulation of the previous integral leads to
$$
 \Delta_\alpha\theta=\int_a^b\frac{1}{\sqrt{z-a}}\frac{1}{\sqrt{b-z}}\frac{1}{\sqrt{1+\varepsilon(z)}}\, dz,\qquad  \varepsilon(z)=\frac{b+a}{z-a}\left[ 1- \frac{\frac{w(b)-w(z)}{b-z}}{\frac{w(b)-w(a)}{b-a}}\right]\ . \label{integdelta}
$$
Denoting by $L=b-a$, $q=\frac{L}{b}$ and applying the change of variable $s=\frac{z-a}{L}$,  the apsidal angle assumes the form \cite{MR3159801} 
\begin{equation}\label{eq:deltatheta}
\Delta_{\alpha}\theta=\int_0^1\frac{1}{\sqrt{s(1-s)}}\frac{1}{\sqrt{1+\mathcal E_{\alpha}(s,q)}}ds
\end{equation}
with
\begin{equation}\label{eps}
{\displaystyle\E_{\al}(s,q)=}\left\{
\begin{array}{ll}
	{\displaystyle 
	\frac{2-q}{qs}\left[ 1-\frac{1-\big(1-q(1-s)\big)^{\alpha}}{(1-s)(1-(1-q)^{\alpha})}\right], }&\alpha\neq 0\\
	\\
	{\displaystyle\frac{2-q}{qs}\left[1-\frac{\log(1-q(1-s))}{(1-s)\log(1-q)} \right],\quad }&\alpha=0.
\end{array}\right.
\end{equation}
 The variable $q=1-\frac{r_-}{r_+}$ shares the properties of the orbital  eccentricity $e=(r_+-r_-)/(r_++r_-)$: it is well defined for all bounded orbits, both open or closed, $q=0$ for circular orbits, $q=1$ for colliding orbits and $q$ is monotonically decreasing in the angular momentum, $ \frac{dq}{d\ell}<0$. Also $\frac{de}{dq}>0$, thus $q$ is a parameter equivalent to $e$.
So, the monotonicity of $\Delta_\al\theta$ with respect to $q$ yields the (opposite) monotonicity  in terms of the angular momentum and the monotonicity with respect $e$ as well. 
Therefore we aim at showing  that $\frac{d\Delta_\al\theta }{dq}<0$ for any $\al\in (-2,1) $ and $q\in (0,1)$.
In proving this statement we restrict ourself to the case $\al\in(0,1)$. Indeed the $\al=0$ case has been already proved in \cite{MR3159801} and the  homogenous potentials  with $\al\in (-2,0)$ are in duality relation to those with $\al\in (0,1)$, \cite{grant}. The duality relates exponent $\al,\bar \al$ satisfying $(2-\al)(2-\bar\al)=4$. In particular it descends that $\Delta_\al\theta=\frac{2-\bar \al}{2}\Delta_{\bar \al}\theta$. Hence, the behaviour of the apsidal angle $\Delta_\al\theta$ with $\al\in(0,1)$ is equivalent to the behaviour  of $\Delta_{\bar \al}\theta$ where $\bar \al\in (-2,0)$ satisfies $(2-\al)(2-\bar\al)=4$. Such a relation extends the identification of the inverse square and linear force problem already known by Newton.

Before proving the monotonicity of the apsidal angle, let us recall  the limiting values and  show the asymptotic values for $q\to 0$.

\section{Asymptotic values of the apsidal angle}

For $\al\in(-2,1)$, the limiting values of the apsidal angle when the eccentricity goes to zero and to 1 are well known, see for instance \cite{1997MNRAS}:
$$
\lim_{q\to 0 }\Delta_\al \theta=\frac{\pi}{\sqrt{2-\al}},\quad\qquad  \lim_{q\to 1 }\Delta_\al \theta=\left\{\begin{array}{ll}
\frac{\pi}{2-\al},\quad &\al\in(0,1)\\
\frac{\pi}{2}, &\al\in(-2,0]\end{array}\right.\ .
$$
These results extend to all power law potentials $V_\al$ with $\al\in (-\infty,2)$.  Assume first that  $\al\in (0,2)$, and compute the limit of \eqref{eq:deltatheta}  for $q\to 0$ and $q\to 1$. Then applying the duality relation we conclude that
$$
\lim_{q\to 0 }\Delta_\al \theta=\frac{\pi}{\sqrt{2-\al}},\ \forall\al\in (-\infty,2),\qquad  \lim_{q\to 1 }\Delta_\al \theta=\left\{\begin{array}{ll}
\frac{\pi}{2-\al},\quad &\al\in(0,2)\\
\frac{\pi}{2}, &\al\in(-\infty,0]\end{array}\right.\ .
$$
Concerning  the asymptotic of $\Delta_\al\theta$ when $q\to 0$, 
for small $q$ we can expand the function $\E_\al(s,q)$ given in \eqref{eps} as $\E_\al(s,q)=(1-\al)+\tilde\E_\al(s,q)$, where 
$$
\begin{aligned}
\tilde\E_\al(s,q)=&\frac{1}{6}(1-\al)(2-\al)\Big[(1-2s)q +\frac{1}{2}(1-s)(\al s-3s +2)q^2\\
&-\frac{1}{60}\left(\begin{array}{l} \left( 6\,{a}^{2}-42\,a+72 \right) {s}^{3}+ \left( -9\,{a}^{2}+93\,a-
198 \right) {s}^{2}\\
+ \left( 182-57\,a+{a}^{2} \right) s+{a}^{2}+3\,a-
58\end{array}\right)q^3\Big]+o(q^3).
\end{aligned}
$$
After inserting into \eqref{eq:deltatheta}, expanding the inverse square root and integrating the series, we obtain the approximation
$$
\Delta_\al\theta=\frac{\pi}{\sqrt{2-\al}}\Big( 1+\frac{1}{96}(\al-1)(\al+2)q^2+\frac{1}{96}(\al-1)(\al+2)q^3 +o(q^3)\Big)
$$
valid for any $\al<2$. Similar to the relations  found in \cite{MNR:MNR8819, MR1433938} in case of $\al\sim 1$, the previous formula shows that for small eccentricity the apsidal angle is decreasing for any $\al\in(-2,1)$.
On the contrary, when $\al<-2$ and $\al\in (1,2)$ the above asymptotic relation implies that the apsidal angle is increasing in the eccentricity. This latter behaviour also agrees with the numerical experiments.

\section{Preliminary results and notation setting} 

In \cite{MR3159801} it has been shown that a sufficient condition for the apsidal angle to be monotonically decreasing with respect to $q$ is that the quantity
\begin{equation}\label{eq:Ial}
\mathcal I_{\al}(s,q)\bydef\partial_{q}\E_{\al}(s,q)\big(1+\E_{\al}(1-s,q)\big)^{2}+\partial_{q}\E_{\al}(1-s,q)\big(1+\E_{\al}(s,q)\big)^{2}
\end{equation}
is positive for any  $s\in (0,\frac{1}{2})$ and $q\in(0,1)$.
Denoting by  
$$
\Omega_{\al,s,q}=(0,1)\times\left(0,\frac{1}{2}\right)\times(0,1).
$$
The ultimate target of the paper is to show the following:
\begin{proposition}\label{prop:main}
For any $(\al,s,q)\in \Omega_{\al,s,q}$ it holds $\mathcal I_\al(s,q)>0$.
\end{proposition}
\begin{remark}\label{rmk:proof}
By proving the Proposition \ref{prop:main} it follows that the apsidal angle is decreasing as a function of $q$ for any $\al\in (0,1)$. Then, the duality condition  discussed above  and $\frac{de}{dq}>0$ provide the proof  of  Theorem \ref{theo:main}.
\end{remark}
Adopting  the  same notation as  in the quoted paper, we introduce the weights 
\begin{equation}\label{eq:w}
\omega^{\alpha}_{n}\bydef
-{\al \choose n}(-1)^{n},  \quad\forall  \al\in(0,1), \qquad \forall n\geq 1.
\end{equation}
For any $n\geq 2,  \al\in(0,1),  s\in(0,1)$, let us define the  functions
\begin{equation}
\begin{aligned}
A_{n}(s)\bydef& \frac{1-(1-s)^{n-1}}{s},\\
K^{\al}_{n}(s)\bydef& 2\frac{\omega^{\al}_{n+1}}{\omega^{\al}_{n}}A_{n+1}(s)-A_{n}(s)=2\frac{n-\al}{n+1}A_{n+1}(s)-A_{n}(s).\label{eq:Kn}
\end{aligned}
\end{equation}
The binomial coefficient ${\al \choose n}$ is equal to $\frac{1}{n !}\al(\al-1)\dots (\al-n+1)$ and extends the usual definition given for integer entries.
Note that $\omega^{\al}_{n}>0$ for any $n\geq 1$ and $\al\in(0,1)$.
\begin{remark}
The functions $K^\al_n(s)$ are continuous at $s=0$. Indeed 
\begin{equation}\label{eq:K_nbis}
K_n^\al(s)=\sum_{k=2}^{n-1}{n-1\choose k}s^{k-1}(-1)^k\left(1-\frac{2(n-\al)}{n+1}\frac{n}{n-k}\right)-2\frac{n-\al}{n+1}s^{n-1}(-1)^n
+\frac{n^2+1-2n\al}{n+1}.
\end{equation}
\end{remark}
By means of the power series expansion
$$
(1-x)^{\al}=1+\sum_{n\geq 1}{\al \choose n}(-1)^{n}x^{n},
$$
the functions $\E_{\al}(s,q)$, $\partial_{q}\E_{\alpha}(s,q)$ can be rewritten in the form
\begin{equation*}%\label{eq:Ea}
\E_{\al}(s,q)=C_{\al}(q)\sum_{n\geq 2}\omega_{n}^{\al}(2-q)q^{n-2}A_{n}(s)
\end{equation*}
\begin{equation}\label{eq:Dedq}
\partial_{q}\E_{\alpha}(s,q)
=C_{\al}(q)^2\sum_{n\geq 2,m\geq 1}\omega^{\al}_{n}\omega^{\al}_{m}\Big((2-q)(n-m)-2\Big)A_{n}(s)q^{(n-2)+(m-2)}
\end{equation}
where
$$
C_{\al}(q)\bydef\left(\sum_{n\geq 1}\omega^{\al}_{n}q^{n-1} \right)^{-1}.
$$
Some results  concerning the properties of $A_n(s)$, $K^\al_n(s)$ and $\E_\al(s,q)$  that will be useful afterwards are now collected. We refer to  \cite{MR3159801} for the proofs.

\begin{lemma}\label{lem:propA}
\begin{itemize}
\item[]
\item[i)] $A_{2}(s)=1$, $A_{n}(1)=1$ for any $n\geq 2$,
\item[ii)] $A_{n+1}(s)-A_{n}(s)=(1-s)^{n-1},\quad A_{n+1}(1-s)-A_{n}(1-s)=s^{(n-1)}$,
\item[iii)] for any $n\geq 3$ and $s\in(0,1)$, $A_{n}(s)$ is decreasing and convex, i.e. $A_{n}'(s)<0$ and $A_{n}''(s)>0$, 
\item[iv)] for any $n>m$, $A'_{n}(s)<A'_{m}(s)$  for any $s\in(0,1)$.
\end{itemize}
\end{lemma}

\begin{lemma}\label{lem:propKn}
For any $\al\in(0,1)$ and $n\geq 2$
\begin{itemize}
\item[i)] $K^{\al}_{n}(s)$ is decreasing and convex  for any $s\in(0,1)$,
\item[ii)] $K^{\al}_{n+1}(s)-K^{\al}_{n}(s)>0$ for any $s\in(0,1)$.
\end{itemize}
For any $\al\in(0,1)$ and $n\geq 4$
\begin{itemize}
\item[iii)] $K^{\al}_{n+1}(s)-K^{\al}_{n}(s)$ is decreasing for any $s\in(0,1)$.
\end{itemize}
\end{lemma}

\begin{lemma}\label{prop:E}
For any $\al\in(0,1)$ and $q\in(0,1)$
\begin{itemize}
\item[i)] $\E_{\alpha}(s,q)> 0$ for any $s\in(0,1)$,
\item[ii)] $\mathcal E_{\alpha}(s,q)$ is monotonically decreasing and convex in $s$ for any $s\in(0,1)$,
\item[iii)] $\partial_{q}\E_{\alpha}\left(\frac{1}{2},q\right)>0$,
\item[iv)] $\partial_q\E_\al(s,q) $ is monotonically  decreasing and convex in $s$ for any $s\in(0,1)$.
\end{itemize} \end{lemma}

Points $iii) $ and $iv)$ of the previous Lemma assure that $\partial_q\E_\al(s,q) $ is positive for any $s\in (0,\frac{1}{2})$ and any $(\al,q)$ in the range of interest. On the contrary, when $s>1/2$ it is not assured that $\partial_q\E_\al(s,q)$ is positive for any values of $\al,q$. In fact, for any $s>1/2$ there is a region of the $\al$-$q$ square where $\partial_q\E_\al(s,q)<0$ and this features makes the proof of the positivity of $\mathcal I_\al(s,q)$ a non trivial problem. The positivity of $\mathcal I_\al(s,q)$ can not be argued from the global behaviour of $\E_{\al}$ and $\partial_{q} \E_{\al}$. Indeed, if for convexity it holds   $\partial_q\E_\al(s,q)>|\partial_q\E_\al(1-s,q)|$,  the same convexity argument implies   $\E_\al(s,q)>\E_\al(1-s,q)$. 
Henceforth,  a  detailed and deep analysis of the function $\mathcal I_\al(s,q)$ is necessary.

Our strategy to prove Proposition \ref{prop:main} combines on-paper analytical calculations  and rigorous computations performed with the assist of a computer.  We proceed as follows: first   we  prove in section \ref{sec:q>09}  that $\mathcal I_\al(s,q)>0$ for any  $\al\in (0,1)$, $s\in (0,1/2)$, and $q\in[0.9,1)$ by showing that $\partial_q\E_\al(s,q)>0$ for any $s$, $ \al$ of interest and $q\in[0.9,1)$. Then, in section \ref{sec:q<09}, we address the case $q<0.9$. The quantity $\mathcal I_\al(s,q)$ will be written as a power series in the variable $q$, i.e. $\mathcal I_\al(s,q)=\sum_p Q_\al^p(s,q)q^p$. The series is then decomposed in the sum of a finite part $\mathcal I^f_\al(s,q)$ ($p\leq \bar p$) and an infinite tail part $\mathcal I^\infty_\al(s,q)$ ($p>\bar p)$. From a certain index $p$ afterwards we show that all the coefficients $Q_\al^p(s,q)$  are positive hence proving the positivity of the infinite tail part. Finally, by rigorous computation,  the finite part is proven to be positive. 

Before proceeding with the proof of the above statements, let us spend some words on what we mean by rigorous computation and to explain the computational technique.

\section{Rigorous computing by means of interval arithmetics}

Very often, when proving a theorem, one produces some numerical computations so to have an evidence of the validity of the statement. However, in most of the cases,  such an evidence can not be awarded as a proof of the theorem. The first simple reason is that any computer is a {\it finite} machine and the   algorithms therein implemented undergo discretisation.  Even the simplest operation, when performed by a computer, may be not exact due to rounding errors. For instance, the standard  floating-point arithmetics would never produce an exact result of the sum $\pi+\sqrt 2$ because both $\pi$ and $\sqrt 2$ can not be exactly represented in the machine. 

One of the goals of a rigorous computational method is to equip the computers with algorithms that ensure  that the output of the computation is reliable and rigorous in the mathematical sense. In the past decades a large number of rigorous computational  algorithms have been produced to address questions  in different areas of mathematics and a variety of  fundamental problems have been solved. Any survey in this direction  is out of the scope of this work and we refer to \cite{MR2652784} for a  review. 

In this paper we will use the computers to prove that certain functions have a determined sign in some domains. For this we combine the floating point arithmetics with the so called {\it Interval Arithmetics}. In regime of interval arithmetics,  numbers are replaced by  intervals and  real number operations are adapted to interval entries.   Given  $\bfA$, $\bfB$ two  intervals and an operation $\circ\in\{ +,-,\cdot,/\}$,  the corresponding operation between intervals, still denoted by $\circ$,  has to satisfy the inclusion property, that is $\bfA\circ \bfB\supset \{ a\circ b,\forall a\in \bfA, \forall b\in \bfB\}$.
Denoting by $\bfA=[\underline a, \bar a]$, $\bfB=[\underline b, \bar b]$ two intervals in $\R$, the sums and the difference are defined as
$$
\bfA+\bfB=[\underline a+\underline b,\bar a +\bar b],\quad \bfA-\bfB=[\underline a-\bar b,\bar a-\underline b].
$$
For the multiplication and division there are different formulas depending on whether $\bfA$, $\bfB$ are positive, negative or contain the zero. Interval extension can  be produced for many of the standard functions, like the exponential, logarithmic,  trigonometric and hyperbolic functions \cite{intervalbook}. 

Suppose $a,b$ are real numbers and we are interested in the result of $a\circ b$.  The basic algorithm  of our rigorous method   is the following: choose intervals $\bfA$, $\bfB$ whose bounds are floating-point numbers such that $a\in \bfA$ and $b\in \bfB$. Taking into account any possible rounding error, define $\bfC$ so that $\bfC\supset \bfA\circ\bfB$. This assures that $a\circ b\in \bfC$. In other words we find rigorous bound of the aimed solution $a\circ b$. Clearly the size of the interval $\bfC$ may increase dramatically, depending on the choice of $\bfA$, $\bfB$, on the nature of the operation and on the number of operations. Nevertheless, the enclosure is rigorous in the mathematical sense, although it could be useless.

Besides resulting in rigorous enclosures, the interval arithmetics is also well  suited for our computations because we  need to study functions defined not just on a finite number of points rather or some bounded domains. Hence, by splitting the domain into the union of a finite number of multi-dimensional  intervals, we can verify the properties of a function on the whole domain with a finite number of computations.

We remark that  any computation in regime of  interval arithmetics results an interval containing the exact solution. Therefore
the method can be adopted to rigorously verify open inequalities on closed domains. In the sequel of the paper we will often encounter problems similar to the following.  Given $f:[0,1]\to \R$  a continuos function  such that   $f(0)=0$, prove that $f(x)>0$ for any $x\in(0,1]$. We would never be able to prove by interval arithmetics that $f(x)>0$ for $x\in(0,1]$. Indeed,  for $x$ small enough, the interval enclosure of $f(x)$ is an interval that has both positive and negative elements. In practice,  the open domain $(0,1]$ is replaced by the union of closed intervals $\cup_i \bfA_i\supset (0,1]$, hence,  whatever domain decomposition we choose, we need to evaluate $f(\bfA)$ with $0\in \bfA$. The result will be an interval $\bfB$ containing zero at least. In this situation, in order to prove that $f(x)>0$,  we first need to factor out powers or functions of $x$, for instance   $f(x)=x^ph(x)$, and prove that $h(x)$ is strictly positive for any $x\in[0,1]$.

\section{Proof of  the monotonicity for almost radial orbits}\label{sec:q>09}
We start with the proof of Proposition \ref{prop:main} by proving that $\mathcal I_{\al}(s,q)>0$ for any $ \al\in(0,1)$, $s\in(0,1/2)$, $q\in[0.9,1)$. It corresponds to the case the orbits are highly eccentric. Since for Lemma  \ref{prop:E},  $\partial_{q}\E_{\al}(s,q)> 0$ for any $(\al,s,q)\in \Omega_{\al,s,q}$,  it remains to prove that $\partial_{q}\E_{\al}(1-s,q)\geq 0$ for any $\al\in(0,1)$, $s\in(0,1/2)$ and $q\in[0.9,1)$. 

\begin{proposition}\label{prop:0.9}
For any $\al\in(0,1)$,  $s\in(0,1/2)$ and $q\in[0.9,1)$, $\partial_{q}\E_{\al}(1-s,q)\geq 0$. 
\end{proposition}
\proof
For Lemma \ref{prop:E} $iv)$, it is enough to show that $\lim_{s\to 1}\partial_{q}\E_{\al}(s,q)>0$ for $q\in[0.9,1)$ and $\al\in(0,1)$. The computation of the derivative of $\E_\al(s,q)$ gives 
$$
\partial_q \E_\al(1,q)=\frac{1}{q(1-q)(1-(1-q)^{\al})^2}N(\al,q)
$$
with
$$
N(\al,q)\bydef\left[-\frac{2}{q}(1-q)\big(1-(1-q)^{\al}\big)^2+\al q(1-q)\big(1-(1-q)^\al\big) +(2-q)\al^2 q(1-q)^{\al} \right].
$$
The task is to prove by rigorous computation  that $N(\al,q)\geq 0$ for any $\al\in(0,1)$ and $q\in[0.9,1)$. 
Since  $N(0,q)=N(1,q)=N(\al,1)=0$, we first need to factor out powers of $\al$, $(1-\al)$ and $(1-q)$. We address the case $\al$ close to zero and $\al$ close to 1 separately.
\paragraph{$\al$ close to 0}
For a choice of $\bar \al\in (0,1)$ it holds
$$
(1-q)^\al>1+\log(1-q)\al+\frac{1}{2}(1-q)^{\bar\al} \log(1-q)^2\al^2,\qquad \forall \al\in(0,\bar\al)
$$
$$
(1-q)^\al<1+\log(1-q)\al+\frac{1}{2} \log(1-q)^2\al^2,\qquad\quad \forall \al\in(0,\bar\al).
$$
By properly  inserting the above inequalities into $N(\al,q)$, one  obtains
$$
\frac{N(\al,q)}{\al^2(1-q)^\al}> N_0(\al,q), \qquad \forall \al\in(0,\bar\al)
$$
where
\begin{equation}\label{eq:N_0}
\begin{aligned}
N_0(\al,q)\bydef& -\frac{2}{q}(1-q)^{1-\al}\log^2(1-q)\left(1+\frac{1}{2}\al\log(1-q) \right)^2\\
&-(1-q)^{1-\al}\log(1-q)\left(1+\frac{1}{2}(1-q)^{\bar\al} \al\log(1-q)\right) +(2-q)q.
\end{aligned}
\end{equation}
The latter function has a criticality at $q=1$, both because $N_0(\al,1)=0$ and because the logarithmic function blows up.
 We solve both problems at the same time. 
We expand the square and collect the term
$$
F(\al,q)\bydef -(1-q)^\frac{1-\al}{4}\log(1-q)
$$
so to rewrite  $N_0(\al,q)$ in the form
$$
\begin{aligned}
N_0(\al,q)=&-\frac{\al^2}{2q}F(\al,q)^4+\frac{2\al}{q}(1-q)^\frac{1-\al}{4}F(\al,q)^3-(1-q)^\frac{1-\al}{2}\left(\frac{2}{q}+\frac{1}{2}\al(1-q)^{\bar\al}\right)F(\al,q)\\
&+(1-q)^\frac{3(1-\al)}{4}F(\al,q)+(2-q)q.
\end{aligned}
$$
In the variable  $q$ the  function $F(\al,q)$ is increasing up to the value $q_{max}=1-e^{-\frac{4}{1-\al}}$ where the maximal value $\frac{4}{e(1-\al)}$ is achieved. For larger $q$ the function  is decreasing and concave.
Thus,  for any $\al\in(0,1)$ and  $q\in [q_{max},1)$
$$
\tilde R(\al,q)\leq F(\al,q) \leq \frac{4}{e(1-\al)}
$$
where $\tilde R(\al,q)$ is the line joining the point $\Big(q_{max},\frac{4}{e(1-\al)}\Big)$ to $(1,0)$.
By inserting these bounds in the previous relation (note that $1+\frac{1}{2}(1-q)^{\bar\al} \log(1-q)\al>0$ for any $q$ and $\al $ in the range we are interested in), 
 it follows 
 $$
 N_0(\al,q)> \tilde N_0(\al,q),\qquad \forall \al\in(0,\bar \al), \ \forall q\in[q_{max},1)
 $$
 where
\begin{equation}\label{eq:tildeN0}
\begin{aligned}
\tilde N_0(\al,q)\bydef& -\frac{\al^2}{2q}\left( \frac{4}{e(1-\al)}\right)^4+\frac{2\al}{q}(1-q)^\frac{1-\al}{4}\tilde R(\al,q)^3\\
&-(1-q)^\frac{1-\al}{2}\left(\frac{2}{q}+\frac{1}{2}\al(1-q)^{\bar\al}\right)\left( \frac{4}{e(1-\al)}\right)^2+(1-q)^\frac{3(1-\al)}{4}\tilde R(\al,q)+(2-q)q.
\end{aligned}
\end{equation}
Hence,  in order to prove that $N(\al,q)$ is nonegative in some strip $\al$ close to zero, we choose the value of $ \al=0.15$. 
Then we prove by  rigorous computation that $N_0(\al,q)>0$ for $(\al,q)\in[0,\bar \al]\times[0.9,q_{max}]$ and  that $\tilde N_0(\al,q)>0$ for $(\al,q)\in [0,\bar \al]\times[q_{max},1]$. We refer to Appendix \ref{sec:AppN0} for more details.

\paragraph{$\al$ close to 1}
 Similarly, when $\al$ is close to 1, chose a value $\hat\al<1$ and consider the bounds   given by the  truncated Taylor expansion pointed  at $\al=1$
 $$
 (1-q)^\al>(1-q)\Big(1-\log(1-q)(1-\al)+\frac{1}{2}\log^2(1-q)(1-\al)^2\Big), \quad \forall \al\in[\hat \al,1)
 $$
  $$
 (1-q)^\al<(1-q)\Big(1-\log(1-q)(1-\al)\Big)+\frac{1}{2}(1-q)^{\hat\al}\log^2(1-q)(1-\al)^2,\quad \forall \al\in[\hat \al,1).
 $$
 These inequalities lead to
 $$
 \frac{N(\al,q)}{1-\al}>(1-q)N_1(\al,q),\qquad \forall \al\in[\hat\al,1)
 $$
 with
 \begin{equation}\label{eq:N1}
 \begin{aligned}
 N_1(\al,q)\bydef &\Big(-2q(1+\al)+\al q^2\Big)+\log(1-q)\Big(-4(1-q)+\al q(1-q)-\al^2 q (2-q)\Big)\\
 &+(1-\al)\log^2(1-q)\left[-\frac{2(1-q)^2}{q}+2(1-q)-\frac{\al}{2} q(1-q)^{\hat\al}+\frac{2-q}{2}q\al^2 \right]\\
 &+(1-\al)^2\left[ \frac{2}{q}(1-q)^2\log^3(1-q)\right]-(1-\al)^3\left[ \frac{(1-q)^2}{2q}\log^4(1-q)\right].
 \end{aligned}
 \end{equation}
 Again we need to solve the indeterminacy at $q=1$ due to $\log(1-q)$.
 Define
 $$
 G(q)\bydef -(1-q)^\frac{\hat\al}{4}\log(1-q),\quad {\rm and}\quad  \hat q\bydef 1-e^{-\frac{4}{\hat\al}}.
 $$
 The function $G(q)$ is positive for $q\in [0,1]$ and  attains its maximum at $q=\hat q$ and $G(\hat q)=\frac{4}{e\hat \al}$. For larger $q$ the function $G(q)$ is decreasing and concave. Therefore
 $$
 R_2(q)\leq G(q)\leq \frac{4}{e\hat \al},\qquad \forall q\in[\hat q,1)
 $$
 where $R_2(q)=-\frac{4}{\hat \al}e^{-1+\frac{4}{\hat \al}}(q-1+e^{-\frac{4}{\hat \al}})+\frac{4}{e\hat \al}$ is the line joining the maximum to the point $(1,0)$.
 Moreover, the function $-\log(1-q)$ is positive and increasing, therefore for any $q\in [\hat q,1]$ it holds $-\log(1-q)\geq -\log(1-\hat q)$. Note that the quantity in the brackets in the second line of \eqref{eq:N1} is positive for any $q>\frac{1}{2}$. Therefore, 
 by means of the previous  bounds, it follows
 $$
 N_1(\al,q)>\tilde N_1(\al,q)\qquad  \forall\al\in(\hat \al,1),\ \forall q\in[\hat q ,1)
 $$
 where
  \begin{equation}\label{eq:N1tilde}
 \begin{aligned}
\tilde  N_1(\al,q)\bydef&\Big(-2q(1+\al)+\al q^2\Big)+R_2(q)(4-\al q)(1-q)^{1-\frac{\hat \al}{4}}\\
 &+(1-\al)R_2(q)^2\left(-\frac{2}{q}(1-q)^{2-\frac{\hat\al}{2}}+2(1-q)^{1-\frac{\hat\al}{2}}+\frac{\al}{2}q(1-q)^\frac{\hat\al}{2}\right)\\
 &-(1-\al)^2\frac{2}{q}(1-q)^{2-\frac{3}{4}\hat\al}\left(\frac{4}{e\hat \al}\right)^3-(1-\al)^3\frac{(1-q)^{2-\hat\al}}{2q}\left(\frac{4}{e\hat \al}\right)^4\\
 &-(2-q)\al q\log(1-\hat q)+(1-\al)\al^2 q\frac{2-q}{2}\log^2(1-\hat q).
  \end{aligned}
 \end{equation}
 Now we are in position to check by rigorous computations  the positivity of $N_1(\al,q)$ for any $(\al,q)\in[\hat\al,1]\times[0.9,1]$. For the choice  $\hat\al=0.8$,  we check that $N_1(\al,q)>0$ for $\al\in(\hat\al,1)\times[0.9,\hat q]$ and that $\tilde N_1(\al,q)>0$ for $(\al,q)\in[\hat\al,1]\times [\hat q,1]$. See Appendix \eqref{sec:AppN1} for details.
 
  \paragraph{The case $\bar\al\leq \al\leq \hat \al $} 
The last computation concerns the function $N(\al,q)$ for $\al\in[\bar\al,\hat\al]$. The only obstruction now stands on the fact that $N(\al,1)=0$. It is enough to factor out $(1-q)^\al$ and prove that 
\begin{equation}\label{eq:Ntilde}
\tilde N(\al,q)\bydef -\frac{2}{q}(1-q)^{1-\al}(1-(1-q)^{\al})^2+\al q(1-q)^{1-\al}(1-(1-q)^\al) +(2-q)\al^2 q
\end{equation}
is positive. 
See  Appendix \eqref{sec:AppN} for the computational details.
\qed

\section{Proof of $\mathcal I_\al(s,q)>0$ for $q<0.9$}\label{sec:q<09}
For $q<0.9$ it is not true that $\partial_q\E_\al(1-s,q)$ is positive for any $\al\in (0,1)$ and $s\in (0,\frac{1}{2})$. In order to prove that $\mathcal I_\al(s,q)$ defined in \eqref{eq:Ial} is positive we need to take into account  the contributions of all the terms.

Formula \eqref{eq:Dedq} leads to the expansion 
\begin{equation}\label{eq:mI}
\mathcal I_{\al}(s,q)=C_{\al}^2(q)\sum_{p\geq 4}Q^p_{\al}(s,q)q^{p-4}
\end{equation}
where 
 \begin{equation*}\label{eq:Qps1}
\begin{split}
Q^p_{\al}(s,q)=&\sum_{\substack{n\geq 2, m\geq 1\\n+m=p}}\omega^{\al}_{n}\omega^{\al}_{m}2(n-m-1)\left[A_{n}(s)\Big(1+\E_{\al}(1-s,q)\Big)^{2}+A_{n}(1-s)\Big(1+\E_{\al}(s,q)\Big)^{2} \right]\\
&-\sum_{\substack{n\geq 2, m\geq 1\\n+m=p-1}}\omega^{\al}_{n}\omega^{\al}_{m}(n-m)\left[A_{n}(s)\Big(1+\E_{\al}(1-s,q)\Big)^{2}+A_{n}(1-s)\Big(1+\E_{\al}(s,q)\Big)^{2} \right].
\end{split}
\end{equation*}
We split \eqref{eq:mI} in the sum of a finite sum and a tail series
\begin{equation}\label{eq:splitting}
\mathcal I_\al(s,q)=C_\al^2(s,q)\left[\mathcal I^f_\al(s,q)+\mathcal I^\infty_\al(s,q) \right]\bydef C_\al^2(s,q)\left[\sum_{p=4}^{11} Q^p_{\al}(s,q)q^{p-4}+\sum_{p> 11}Q^p_{\al}(s,q)q^{p-4}\right].
\end{equation}
In the following both $\mathcal I^f_\al(s,q)$ and $\mathcal I^\infty_\al(s,q)$ are proven to be positive for any $\al\in(0,1), s\in (0,\frac{1}{2})$, $q\in(0,0.9]$, the latter one by means of analytical estimates, the former by rigorous computations.

\subsection{Analytical bound for the infinite tail series.}
The target of this  section is to show that the coefficients $Q^p_\al(s,q)$ are all positive for $p$ large enough for any $(\al,s,q)\in\Omega_{\al,s,q}$. Let us first rewrite $Q^p_\al(s,q)$ in a more convenient form
\begin{equation*}
\begin{split}
Q^p_{\alpha}(s,q)=&\sum_{\substack{n\geq 1, m\geq 1\\n+m=p-1}}\w^\al_{n+1}\w^\al_{m}2(n-m)\left(A_{n+1}(s)\Big(1+\E_\al(1-s,q)\Big)^{2}+A_{n+1}(1-s)\Big(1+\E_\al(s,q)\Big)^{2} \right)\\
&-\sum_{\substack{n\geq 2, m\geq 1\\n+m=p-1}}\w^\al_{n}\w^\al_{m}(n-m)\left(A_{n}(s)\Big(1+\E_\al(1-s,q)\Big)^{2}+A_{n}(1-s)\Big(1+\E_\al(s,q)\Big)^{2} \right)\\
=&\ \w^\al_{2}\w^\al_{p-2}2(1-(p-2))\left( A_{2}(s)\Big(1+\E_\al(1-s,q)\Big)^{2}+A_{2}(1-s)\Big(1+\E_\al(s,q)\Big)^{2} \right)\\
&+\sum_{\substack{n\geq 2, m\geq 1\\n+m=p-1}}\w^\al_{n+1}\w^\al_{m}2(n-m)\left(A_{n+1}(s)\Big(1+\E_\al(1-s,q)\Big)^{2}+A_{n+1}(1-s)\Big(1+\E_\al(s,q)\Big)^{2} \right)\\
&-\sum_{\substack{n\geq 2, m\geq 1\\n+m=p-1}}\w^\al_{n}\w^\al_{m}(n-m)\left(A_{n}(s)\Big(1+\E_\al(1-s,q)\Big)^{2}+A_{n}(1-s)\Big(1+\E_\al(s,q)\Big)^{2} \right).
\end{split}
\end{equation*}
Collecting all the contributions and recalling the definition \eqref{eq:Kn}
it results
\begin{equation}\label{eq:Qps}
\begin{split}
Q^p_{\al}(s,q)=&\sum_{\substack{n\geq 2, m\geq 1\\n+m=p-1}}\w^\al_{n}\w^\al_{m}(n-m)\left[
K^\al_{n}(s)\Big(1+\E_\al(1-s,q)\Big)^{2}+K^\al_{n}(1-s)\Big(1+\E_\al(s,q)\Big)^{2} \right]\\
&-2\w^\al_2\w^\al_{p-2}(p-3)\Bigg( \Big(1+\E_\al(1-s,q)\Big)^{2} +\Big(1+\E_\al(s,q)\Big)^{2} \Bigg)\\
=&\ T^p(\al,s)\Big(1+\E_\al(1-s,q)\Big)^{2}+T^p(\al,1-s)\Big(1+\E_\al(s,q)\Big)^{2}
\end{split}
\end{equation}
where 
\begin{equation}\label{eq:Tps}
T^p(\al,s):=\sum_{\substack{n\geq 2, m\geq 1\\n+m=p-1}}\w^\al_{n}\w^\al_{m}(n-m)
K^\al_{n}(s)-2\w^\al_2\w^\al_{p-2}(p-3).
\end{equation}
The positivity of $Q^p_\al(s,q)$ in its domain of definition is guaranteed by proving that $T^p(\al,s)$ is positive for any $s\in(0,1)$ and $\al\in (0,1)$.  For that,  
further properties of $K^\al_n(s)$, meant to complete the ones collected in Lemma \ref{lem:propKn}, are necessary.
\begin{lemma}\label{lemma:extra}
For any $\al\in(0,1)$ and $s\in(0,1)$
\begin{itemize}
\item[i)] $K^\al_5(s)-K^\al_3(s)\geq K^\al_5(1)-K^\al_3(1)$,
\item[ii)] for any $n\geq 3$ the function
$$
f(s)\bydef\w^\al_n\w^\al_1(n-1)K^\al_n(s)-(n-3)\w^\al_{n-1}\w^\al_2K^\al_2(s)
$$
is decreasing.
\end{itemize}
\end{lemma} 
\proof
$i)$ Straightforward computations give 
$$
\frac{\Big( K^\al_5(s)-K^\al_3(s)\Big)-\Big( K^\al_5(1)-K^\al_3(1)\Big)}{1-s}=\frac{ 2-s}{6}  \left( 10\,{s}^{2}-2\,a{s}^{2}+4\,as-14\,s
+5-a \right) 
$$
that is positive for any $s\in(0,1)$ and $\al\in (0,1)$.

$ii)$
Up to a positive constant
$$
f'(s)=\frac{n-1-\al}{n}(n-1)(K^\al_n(s))'+\frac{1}{3}(n-3)(1-\al)(2-\al).
$$
The convexity of $K^\al_n(s)$ and the definition of $A_n(s)$ imply that
$$
f'_n(s)<-\frac{n-1-\al}{n}(n-1)\frac{n-2\al-1}{n+1}+\frac{1}{3}(n-3)(1-\al)(2-\al).
$$
By direct computation we prove that $f_n'(s)<0$ for any $\al\in(0,1)$, $s\in(0,1)$ and $n=3,\dots,10$. The estimate
\begin{equation*}
\begin{split}
f'_n(s)&<\frac{1}{3n(n+1)}\Big(-3(n-2)(n-1)(n-3)+(n-3)2n(n+1) \Big)\\
&<\frac{n-3}{3n(n+1)}(-n^2+11n-6)
\end{split}
\end{equation*}
gives $f'_n(s)<0$ for any $n\geq 11$.
\qed

We can now prove that $T^p(\al,s)$ is positive for $p$ large enough.  
\begin{lemma}\label{lem:bound1}
For any $p\geq 4$
\begin{equation}\label{eq:lemma}
\sum_{\substack{n\geq 2, m\geq 1\\n+m=p-1}}\w^\al_{n}\w^\al_{m}(n-m)K^{\al}_{n}(s)\geq\sum_{\substack{n\geq 2, m\geq 1\\n+m=p-1}}\w^\al_{n}\w^\al_{m}(n-m)\frac{n-1-2\al}{n+1},\quad \forall s\in(0,1).
\end{equation}
\end{lemma}
\proof
Set
$$
g_p(s)\bydef \sum_{\substack{n\geq 2, m\geq 1\\n+m=p-1}}\w^\al_{n}\w^\al_{m}(n-m)K^{\al}_{n}(s).
$$
The right hand side of \eqref{eq:lemma} is nothing but the evaluation of the left side at $s=1$, hence  we aim at proving that $g_p(s)\geq g_p(1)$. 
In the cases $p=4$ and $p=5$ all the quantities $(n-m)$ appearing in the series are nonnegative, hence the result follows directly from the monotonicity of $K^\al_n(s)$.

 For $p=6$,  $g_6(s)=\w^\al_3\w^\al_2(K^\al_3(s)-K^\al_2(s))+\w^\al_4\w^\al_1 K^\al_4(s)$  consists in a cubic polynomial in $s$ with $\al$-dependent  coefficients.  Up to positive factors, one computes $g_6'''(s)= -(\al^2-7\al+12)$ showing that $g_6'''(s)<0$ for any $s\in(0,1)$ and any $\al\in (0,1)$. Then the straightforward computation of $g_6''(1)>0$ and $g_6'(1)<0$ allows to conclude that $g_6(s)$ is monotonically decreasing.
 The same argument applies for $p=8$. In that case $g_8(s)$ is a quintic polynomial then one starts computing from the  fifth derivative.
 
For the values  $p=7$ and $p\geq 9$ we proceed as  follows. Assume $p$ is odd, $p=2\ell+1$.
Then
\begin{equation*}
\begin{split}
g(s)=&\sum_{\substack{n=\ell+1\\(m=2\ell-n)}}^{2\ell-3}\w^\al_n\w^\al_m(n-m)\Big(K^{\al}_{n}(s)-K^{\al}_{m}(s) \Big)\\
&+\w^\al_{2\ell-2}\w^\al_{2}(2\ell-4)\Big(K^{\al}_{2\ell-2}(s)-K^{\al}_{2}(s)\Big)+\w^\al_{2\ell-1}\w^\al_1(2\ell-2)K^{\al}_{2\ell-1}(s).\\
\end{split}
\end{equation*}
Any term $(K^{\al}_{n}(s)-K^{\al}_{m}(s))$ inside the series is such that $m\geq 3$ and $n\geq m+2$.
We can write
\begin{equation*}
\begin{array}{ll}
{\displaystyle K^\al_n(s)-K^\al_m(s)=\sum_{i=m}^{n-1}K^\al_{i+1}(s)-K^\al_i}(s)& {\rm if}\ m>3,\\
{\displaystyle K^\al_n(s)-K^\al_m(s)=\sum_{i=5}^{n-1}K^\al_{i+1}(s)-K^\al_i(s)+K^{\al}_{5}(s)-K^{\al}_{3}(s)}& {\rm if }\ m=3.
\end{array}
\end{equation*}
In both the cases, for point $iii)$  of Lemma \ref{lem:propKn} and  Lemma \ref{lemma:extra}$i)$ we infer
$$
K^{\al}_{n}(s)-K^{\al}_{m}(s)\geq K^{\al}_{n}(1)-K^{\al}_{m}(1).
$$
Point $i)$ of Lemma \ref{lem:propKn} implies
$
K^{\al}_{2\ell-2}(s)\geq K^{\al}_{2\ell-2}(1)
$
 and, from Lemma \ref{lemma:extra}$ii)$, 
 $$
 \w^\al_{2\ell-1}\w^\al_1(2\ell-2)K^{\al}_{2\ell-1}(s)-\w^\al_{2\ell-2}\w^\al_{2}(2\ell-4)K^\al_2(s)\geq  \w^\al_{2\ell-1}\w^\al_1(2\ell-2)K^{\al}_{2\ell-1}(1)-\w^\al_{2\ell-2}\w^\al_{2}(2\ell-4)K^\al_2(1)
$$ for any $\ell\geq 1$. We then conclude that $g(s)\geq g(1)$. 
The case $p$ even is completely equivalent.
\qed

\begin{lemma}\label{lem:bound2}
For any $p\geq 11$ 
$$
\sum_{\substack{n\geq 2, m\geq 1\\n+m=p-1}}\w^\al_{n}\w^\al_{m}(n-m)K^{\al}_{n}(s)\geq 2\w^\al_2\w^\al_{p-2}(p-3), \quad \forall \al\in(0,1), \ s\in (0,1).
$$
\end{lemma}
\proof
From the previous Lemma
$$
\begin{aligned}
\sum_{\substack{n\geq 2, m\geq 1\\n+m=p-1}}\w^\al_{n}\w^\al_{m}(n-m)K^{\al}_{n}(s)&\geq\sum_{\substack{n\geq 2, m\geq 1\\n+m=p-1}}\w^\al_{n}\w^\al_{m}(n-m)\frac{n-1-2\al}{n+1}\\
&\geq\sum_{\substack{n\geq 2, m\geq 1\\n+m=p-1}}\w^\al_{n}\w^\al_{m}(n-m)-\sum_{\substack{n\geq 2, m\geq 1\\n+m=p-1}}\w^\al_{n}\w^\al_{m}(n-m)\frac{2(1+\al)}{n+1}.
\end{aligned}
$$
For symmetry reasons, in the first sum all but the $(m=1)-$term vanish, hence it remains 
\begin{equation}\label{eq:maior}
\sum_{\substack{n\geq 2, m\geq 1\\n+m=p-1}}\w^\al_{n}\w^\al_{m}(n-m)=\w^\al_{p-2}w^\al_1(p-3).
\end{equation}
The second sum amounts to
\begin{equation*}
\begin{aligned}
\sum_{\substack{n\geq 2, m\geq 1\\n+m=p-1}}\w^\al_{n}&\w^\al_{m}(n-m)\frac{2(1+\al)}{n+1}=\frac{2(1+\al)}{p-1}(p-3)\w^\al_{p-2}\w^\al_1+\sum_{\substack{n\geq 2, m\geq 2\\n+m=p-1}}\w^\al_{n}\w^\al_{m}(n-m)\frac{2(1+\al)}{n+1}.
\end{aligned}
\end{equation*}
Since 
$$
\begin{aligned}
\sum_{\substack{n\geq 2, m\geq 2\\n+m=p-1}}\w^\al_{n}\w^\al_{m}\frac{n-m}{n+1}&=\sum_{\substack{n>m\geq 2\\n+m=p-1}}\w^\al_{n}\w^\al_{m}\frac{n-m}{n+1}+\sum_{\substack{m>n\geq  2\\n+m=p-1}}\w^\al_{n}\w^\al_{m}\frac{n-m}{n+1}\\
&=\sum_{\substack{n>m\geq 2\\n+m=p-1}}\w^\al_{n}\w^\al_{m}\frac{n-m}{n+1}+\sum_{\substack{n>m\geq 2\\n+m=p-1}}\w^\al_{m}\w^\al_{n}\frac{m-n}{m+1}\\
&=\sum_{\substack{n>m\geq 2\\n+m=p-1}}\w^\al_{n}\w^\al_{m}(n-m)\left(\frac{1}{n+1}-\frac{1}{m+1}\right),
\end{aligned}
$$
it follows 
\begin{equation}\label{eq:minor}
\begin{aligned}
\sum_{\substack{n\geq 2, m\geq 1\\n+m=p-1}}\w^\al_{n}\w^\al_{m}(n-m)\frac{2(1+\al)}{n+1}=&\frac{2(1+\al)}{p-1}(p-3)\w^\al_{p-2}\w^\al_1\\
&+2(1+\alpha)\sum_{\substack{n>m\geq 2\\n+m=p-1}}\w^\al_{n}\w^\al_{m}(n-m)\left(\frac{1}{n+1}-\frac{1}{m+1}\right).\\
\end{aligned}
\end{equation}

Combining  \eqref{eq:maior} and \eqref{eq:minor}, the thesis follows by proving that
$$
2(1+\al)\sum_{\substack{n>m\geq 2\\n+m=p-1}}\w^\al_{n}\w^\al_{m}(n-m)\left(\frac{1}{m+1}-\frac{1}{n+1}\right)+\w^\al_{p-2}(p-3)\al\left(\al-\frac{2(1+\al)}{p-1}\right)\geq 0.
$$
Note that any contribution of the   sum is positive for  any $\al\in(0,1)$. However  for any $p\geq 4$ there is an interval of $\al$ where the  second term is negative. Therefore the idea is to extract  few elements from the  sum  enough to compensate the negative contribution of the second term.

Assume $p\geq 11$ and rewrite the  sum as
$$
\sum_{\substack{m=2\\n+m=p-1}}^4\w^\al_{n}\w^\al_{m}(n-m)\left(\frac{1}{m+1}-\frac{1}{n+1}\right)+\sum_{\substack{n>m\geq 5\\n+m=p-1}}\w^\al_{n}\w^\al_{m}(n-m)\left(\frac{1}{m+1}-\frac{1}{n+1}\right).
$$
It is enough to show that 
$$
F(\al,p)\bydef 2(1+\al)\sum_{\substack{m=2\\n+m=p-1}}^4\w^\al_{n}\w^\al_{m}(n-m)\left(\frac{1}{m+1}-\frac{1}{n+1}\right)-\w^\al_{p-2}(p-3)\al\left(\frac{2(1+\al)}{p-1}-\al\right)\geq 0
$$
for any $\al\in (0,1)$ and $p\geq 11$.
Writing explicitly the elements, 
$$
\begin{aligned}
F(\al,p)=&2(1+\al)\left[\w^\al_{p-3}\w^\al_2\frac{(p-5)^2}{3(p-2)} +\w^\al_{p-4}\w^\al_3\frac{(p-7)^2}{4(p-3)}+\w^\al_{p-5}\w^\al_4\frac{(p-9)^2}{5(p-4)}\right]\\
&-\w^\al_{p-2}(p-3)\al\left(\frac{2(1+\al)}{p-1}-\al\right).
\end{aligned}
$$
Applying recursively  the relation $\frac{\w^\al_n}{\w^\al_{n-1}}=\frac{n-1-\al}{n}$,  it follows
$$
\frac{\w^\al_{p-4}}{\w^\al_{p-5}}=\frac{(p-5-\al)}{(p-4)},\ \frac{\w^\al_{p-3}}{\w^\al_{p-5}}=\frac{(p-4-\al)(p-5-\al)}{(p-3)(p-4)},\  \frac{\w^\al_{p-2}}{\w^\al_{p-5}}=\frac{(p-3-\al)(p-4-\al)(p-5-\al)}{(p-2)(p-3)(p-4)}.
$$
We can then  factor out the term $w^\al_{p-5}$ and by straightforward computations we replace $F(\al,p)$ with the following polynomial ( still denoted by $F(\al,p)$)
$$
F(\al,p)=\begin{array}{l}
-884\,{p}^{4}+7860\,{p}^{3}-31700\,{p}^{2}+58344\,p-39416+36\,{p}^{5}\\+
 \left( -49920-28620\,{p}^{2}+60450\,p+6940\,{p}^{3}-900\,{p}^{4}+50\,
{p}^{5} \right) \al\\+ \left( 7855\,{p}^{2}-1000\,p-12380-3405\,{p}^{3}+
565\,{p}^{4}-35\,{p}^{5} \right) {\al}^{2}\\+ \left( 10\,{p}^{5}-180\,{p}^
{4}+1280\,{p}^{3}-4140\,{p}^{2}+6870\,p-6240 \right) {\al}^{3}\\+ \left( -
135\,{p}^{3}+385\,{p}^{2}-464\,p-44+19\,{p}^{4}-{p}^{5} \right) {\al}^{4
}.
\end{array}
$$
We start proving that the fourth derivative of $F(\al,p)$ with respect to $\alpha$ is negative for any $\al$ and  any $p\geq 11$. Up to positive factors,  such a derivative amounts  in the polynomial $f(p)= -135\,{p}^{3}+385\,{p}^{2}-464\,p-44+19\,{p}^{4}-{p}^{5} $. Since $f^{(v)}(p)<0$ for any $p$ and  the first up to the fourth order derivatives of $f(p)$ are negative when evaluated at $p=11$, it follows that $f(p)<0$ for any $p\geq 11$.

The lower order $\alpha$-derivative of $F(\al,p)$ evaluated in $\al=1$ are 
$$
\frac{d^3 F}{d\al^3}(1,p)=36p^5-624p^4+4440p^3-15600p^2+30084p-38496,
$$
$$
\frac{d^2 F}{d\al^2}(1,p)=-4510p^2+33652p-62728-750p^3+278p^4-22p^5.
$$
Arguing as before, it is easy to prove that, for any $p\geq 11$,  $\frac{d^3 F}{d\al^3}(1,p)>0$, $\frac{d^2 F}{d\al^2}(1,p)<0$. From where it follows that $\frac{d^2 F}{d\al^2}(\al,p)<0$ for any $\al\in(0,1)$ and $p\geq 11$. Unfortunately it is not true that $\frac{d F}{d\al}(1,p)>0$ for any $p\geq 11$. In that case   the function $F(\al,p)$   would be $\al$-monotone  and it would be enough to check $F(0,p)>0$. However, since $F(\al,p)$ is concave, the positivity of $F(\al,p)$ for any $\al\in(0,1)$ and $p\geq 11$ follows by proving that $F(1,p)-F(0,p)>0$ for any $p\geq 11$.

It holds
$$
F(1,p)-F(0,p)=-68584-24520p^2+65856p+4680p^3-496p^4+24p^5
$$
The fifth derivative of the right hand side  is positive and all the lower-order derivates are positive at $p=11$. Hence $F(1,p)-F(0,p)>0$ for any $p\geq 11$. In conclusion $F(\al,p)>0$ for any $\al\in(0,1)$ and $p\geq 11$.
\qed

\begin{proposition}\label{theo:infinite}
The series $\mathcal I^\infty_\al(s,q)$  is positive for any $(\al,s,q)\in \Omega_{\al,s,q}$.
\end{proposition}
\proof
From the previous Lemma it descends that 
$$
Q^p_\al(s,q)>0,\qquad  \forall (\al,s,q)\in\Omega_{\al,s,q}
$$
for any $p\geq 11$.
Therefore $\mathcal I^\infty_\al(s,q)>0$ for all $(\al,s,q)\in \Omega_{\al,s,q}$.
\qed

We have proved  that the infinite tail series $\mathcal I^\infty_\al(s,q)$ in the splitting \eqref{eq:splitting} is positive for any $(\al,s,q)\in \Omega_{\al,s,q}$. In the next section we are concerned with the finite sum $\mathcal I^f_\al(s,q)$.

\subsection{Rigorous computation of  the finite part}

In order to prove that $\mathcal I_\al(s,q)>0$ for any  $(\al,s,q)\in\Omega_{\al,s,q}$ it remains to prove that 
$$
\mathcal I_\al^f(s,q)\bydef\sum_{p=4}^{10}Q^p_\al(s,q)q^{p-4}
$$
is positive for any $\al\in(0,1)$, $s\in(0,\frac{1}{2})$ and $q\in(0,0.9)$.
Denote by 
$$\Lambda_{\al,s,q}=(0,1)\times\left(0,\frac{1}{2}\right)\times(0,0.9),\quad {\rm and}\quad  \bar \Lambda=[0,1]\times\left[0,\frac{1}{2}\right]\times[0,0.9].$$
Recalling \eqref{eq:Qps}, we can write
$$
\mathcal I_\al^f(s,q)=R_\al(s,q)(1+\E_\al(1-s,q))^2+R_\al(1-s,q)(1+\E_\al(s,q))^2
$$
where
$$
R_\al(s,q)\bydef\sum_{p=4}^{10}\left[ \sum_{\substack{n\geq 2, m\geq 1\\n+m=p-1}}\w^\al_{n}\w^\al_{m}(n-m)K^\al_n(s)-2\w^\al_2\w^\al_{p-2}(p-3)\right]q^{p-4}.
$$
Note that  in the square bracket of the previous relation we have nothing else but the coefficient $T^p(\al,s)$ \eqref{eq:Tps}.

Let us rephrase  what we did and what we plan to do. Previously we proved that $T^p(\al,s)>0$ for any $\al\in(0,1)$ and $s\in (0,1)$ whenever $p\geq 11$. If we were able to prove that $T^p(\al,s)>0$ for any $p\geq 4$ and $s\in(0,1)$, the proof of Proposition \ref{prop:main} would be done. We knew since the beginning  that this could not be the case, since we are aware of the fact  that the derivative $\partial_q\E_\al(s,q)$ is negative for some values $(\al,q)$  when $s>\frac{1}{2}$. Anyway,  the analysis of the previous section is illuminating because it provides the  separation of  the possible negative contributions of $\partial_q\E_\al(s,q)$ from the positive ones. Moreover, a careful analysis of the values of $K^\al_n(\frac{1}{2})$ together with the monotonicity of $K^\al_n(s)$  allows to conclude that for $s\in(0,\frac{1}{2})$ the terms  $T^p(\al,s)$ are all positive, even for $p<11$.

This means that for any $s\in(0,\frac{1}{2})$,  $R_\al(s,q)$ is positive  while $R_\al(1-s,q)$ is negative for some values of $\al,q$. Henceforth the idea is to check that, whereas $R_\al(1-s,q)$ is negative, the positive contribution of $R_\al(s,q)$ is large enough to compensate the negative contribution due to   $R_\al(1-s,q)$.
Since $\mathcal I^f_\al(s,q)$ is a finite sum and its evaluation  amounts to a finite number of algebraic computations, we prove the positivity of $\mathcal I_\al^f(s,q)$ by rigorous computation.

As before, a preliminary analysis  on the quantity $\mathcal I ^f_\al(s,q)$ is necessary.
 Indeed,  it is evident  that $R_\al(s,q)$ and hence $\mathcal I^f_\al(s,q)$ is zero at $\al=1$ and $\al=0$. Moreover, it will be clear later  that $\mathcal I_\al(s,q)$ tends to zero as $q\to 0$. Therefore we first need to factor out powers of $\al$, $(1-\al)$ and $q$ and then to prove that the remaining sum is strictly positive.
 
 This being said, for any $n\geq 2$  introduce $\tilde\w^\al_n$ so that
 $$
 \w^\al_n=\al(1-\al)\tilde\w^\al_n.
 $$
Note that $\tilde\w^\al_n>0$ for any $n\geq 2$ and $\al\in[0,1]$. Rewrite $R_\al(s,q)$ as follows
\begin{equation*}
\begin{split}
R_\al &(s,q)=\sum_{p=4}^{10}\Bigg[ \sum_{\substack{n\geq 2, m\geq 2\\n+m=p-1}}\w^\al_{n}\w^\al_{m}(n-m)K^\al_n(s)+\w^\al_{p-2}\al(p-3)K^\al_{p-2}(s)-\al(1-\al)\w^\al_{p-2}(p-3)\Bigg]q^{p-4}\\
&=\al^2(1-\al)\sum_{p=4}^{10}\left[ \sum_{\substack{n\geq 2, m\geq 2\\n+m=p-1}}(1-\al)\tilde\w^\al_{n}\tilde\w^\al_{m}(n-m)K^\al_n(s)\right.\\
&\qquad \qquad\qquad\qquad   +\tilde\w^\al_{p-2}(p-3)K^\al_{p-2}(s)-(1-\al)\tilde\w^\al_{p-2}(p-3)\Bigg]q^{p-4}\\
&=\al^2(1-\al)\sum_{p=4}^{10}\left[ \sum_{\substack{n\geq 2, m\geq 2\\n+m=p-1}}(1-\al)\tilde\w^\al_{n}\tilde\w^\al_{m}(n-m)K^\al_n(s) +\tilde\w^\al_{p-2}(p-3)\Big(K^\al_{p-2}(s)-(1-\al)\Big)\right]q^{p-4}.\\
\end{split}
\end{equation*}
Since $K^\al_n(s)$ is increasing in $n$, see Lemma \ref{lem:propKn}, the innermost sum is positive for any $p$. Therefore  we have that
$$
R_\al(s,q)\geq \al^2(1-\al)\tilde R_\al(s,q),\quad \forall (\al,s,q)\in(0,1)\times(0,1)\times (0,1) 
$$
where
$$
\widetilde R_\al(s,q)\bydef\sum_{p=4}^{10}\Big[ \tilde\w^\al_{p-2}(p-3)\Big(K^\al_{p-2}(s)-(1-\al)\Big)\Big]q^{p-4}.
$$
Hence, in order to conclude that $\mathcal I^f_\al(s,q)$ is positive for any $(\al,s,q)\in\Lambda_{\al,s,q}$, it is enough to prove that
$$
\widehat {\mathcal I}_\al^f(s,q)\bydef \widetilde R_\al(s,q)(1+\E_\al(1-s,q))^2+\widetilde R_\al(1-s,q)(1+\E_\al(s,q))^2
$$
is positive in $\Lambda_{\al,s,q}$. Note that  $\widetilde R_\al(s,q)$ has less terms  than $R_\al(s,q)$ hence  reducing the number of interval computation and, more important, $\widehat {\mathcal I}_\al^f(s,q)$ does not tend to zero as $\al\to 1$ and $\al\to 0$,  $\tilde\w^0_n$ and $\tilde\w^1_n$ being non zero. 

To further reduce the number of interval computations and, as a consequence, to obtain a better enclosure of the evaluation of $\widehat {\mathcal I}^f_\al(s,q)$, 
we remove the dependence of $s$ in the functions $\E_\al(1-s,q)$, $\E_\al(s,q)$.

Since $\widetilde R_\alpha(1-s,q)$ is supposed to be  negative, the monotonicity of $\E_\al(s,q)$ (see Lemma \ref{prop:E}) implies that
$$
\widehat {\mathcal I}^f_\al(s,q)\geq \widetilde {\mathcal I}^f_\al(s,q)
$$
\begin{equation}\label{eq:tildeIfasq}
\widetilde {\mathcal I}^f_\al(s,q)\bydef \widetilde R_\al(s,q)(1+\E_\al(1,q))^2+\widetilde R_\al(1-s,q)(1+\E_\al(0,q))^2.
\end{equation}
The functions $\E_\al(1,q)$ and $\E_\al(0,q)$ are given below:
$$
\E_\al(0,q)=\frac{ \left( 2-q \right) }{q} \left(-1+\frac{ \al q \left( 1-q \right) ^{-1+\al}}{1-
 \left( 1-q \right) ^{\al}} \right) 
$$
$$
\E_\al(1,q)=\frac{(2-q)}{q}\left(1-\frac{\al q}{1-(1-q)^\al}\right).
$$
The problem of checking the positivity of $\mathcal I^f_a(s,q)$ reduces then to the problem of checking that $\widetilde {\mathcal I}^f_\al(s,q)>0$ for any $(\al,s,q)\in \Lambda_{\al,s,q}$. We still have some problems to solve before doing the rigorous computation. The first one concerns the $q$ variable. Indeed
$$
\widetilde R_\al(0,s)=\frac{1}{6}(1-2s)(2-a)=-\widetilde R_\al(0,1-s)
$$
and both 
$$
\E_\al(0,q), \E_\al(1,q)\to 1-\al, \quad {\rm as}\ q\to 0
$$
meaning that $\widetilde{\mathcal I}^f_\al(s,q)\to 0$ for any $\al\in(0,1)$, $s\in(0,\frac{1}{2})$ as $q\to 0$.
Moreover each of the terms $\E_\al(0,q)$ and $\E_\al(1,q)$ has  $q$ in the denominator, hence an indeterminacy arises when $q\to 0$. 
The second issue consists in the indeterminacy of $\E_\al(0,q)$ and $\E_\al(1,q)$ when $\al=0$.

 To deal with these problems, for a choice of $\bar q$ and $\alpha_0$, we subdivide the domain so that  $\bar\Lambda_{\al,q,s}=\Lambda^0_{\al,s,q}\cup \Lambda^1_{\al,s,q}\cup\Lambda^3_{\al,s,q}$ where
 $$
 \Lambda^0_{\al,s,q}=\left\{(\al,s,q)\in[ 0,1]\times\left[0,\frac{1}{2}\right]\times[0,\bar q]\right\},\quad  \Lambda^1_{\al,s,q}=\left\{(\al,s,q)\in[ 0,\al_0]\times\left[0,\frac{1}{2}\right]\times[\bar q,0.9]\right\}.
 $$
To overcome the indeterminacy due to $q=0$ and $\al=0$,  two different bounds for  $\widetilde{\mathcal I}^f_\al(s,q)$ are found on the  domains $\Lambda^0_{\al,s,q}$ and  $\Lambda^1_{\al,s,q}$.

 \subsubsection{Bound for $\widetilde {\mathcal I}^f_\al(s,q)$ in $\Lambda^0_{\al,s,q}$ }
 
 A truncated Taylor series at $q=0$ provides the estimates 
 $$
 \E_\al(0,q)<(1-\al)+q E_0(q),\qquad \E_\al(1,q)>(1-\al)+q E_1(q),\ \forall q\in[0,\bar q]
$$
with
$$
E_0(q)\bydef \frac{1}{6}(1-\al)(2-\al)\Big(1+(1-\bar q)^{\al-3}q \Big)
$$
$$
E_1(q)\bydef -\frac{1}{6}(1-\al)(2-\al).
$$
Then, for those $s$ where $\widetilde R_\al(1-s,q)<0$ and $ \forall \al\in(0,1), q\in(0,\bar q)$, it holds 
$$
\begin{aligned}
\widetilde {\mathcal I}_\al^f(s,q)&>\widetilde R_\al(s,q)\Big(1+(1-\al)+ qE_1(q)\Big)^2 +\widetilde R_\al(1-s,q)\Big(1+(1-\al)+qE_0(q)\Big)^2 \\
&\geq \frac{1}{6}(1-2s)(2-\al)\Big[\Big(2-\al+q E_1(q)\Big)^2-\Big(2-\al+q E_0(q)\Big)^2\Big]\\
&\ \ +\sum_{p=5}^{10}\Big[ \tilde\w^\al_{p-2}(p-3)(K^\al_{p-2}(s)-(1-\al))\Big]q^{p-4}\Big((2-\al)+qE_1(q)\Big)^2\\
&\ \ +\sum_{p=5}^{10}\Big[ \tilde\w^\al_{p-2}(p-3)(K^\al_{p-2}(1-s)-(1-\al))\Big]q^{p-4}\Big((2-\al)+qE_0(q)\Big)^2. \
\end{aligned}
$$
The $0$-th order term vanishes, and dividing by $q$ we have
\begin{equation*}%\label{eq:Maf}
\begin{aligned}
\frac{\widetilde {\mathcal I}_\al^f(s,q)}{q}>M^f_\al(s,q)\bydef & \frac{1}{6}(1-2s)(2-\al)\left[2(2-\al)\Big( E_1(q)-E_0(q)\Big)+q\Big( E_1^2(q)-E_0^2(q)\Big)\right]\\
&\  +\sum_{p=5}^{10}\Big[ \tilde\w^\al_{p-2}(p-3)(K^\al_{p-2}(s)-(1-\al))\Big]q^{p-5}\Big((2-\al)+qE_1(q)\Big)^2\\
&\  +\sum_{p=5}^{10}\Big[ \tilde\w^\al_{p-2}(p-3)(K^\al_{p-2}(1-s)-(1-\al))\Big]q^{p-5}\Big((2-\al)+q E_0(q)\Big)^2.\\
\end{aligned}
\end{equation*}
Hence,  the positivity of $
M_\al^f(s,q)$ in $\Lambda^0_{\al,s,q}$ implies the positivity of  $\widetilde{\mathcal I}^f_\al(s,q)$ and, as a consequence, of  $\mathcal I^f_\al(s,q)$ for $(\al,s,q)\in(0,1)\times(0,\frac{1}{2})\times(0,\bar q]$.

 \subsubsection{Bound for $\widetilde {\mathcal I}^f_\al(s,q)$ in $\Lambda^1_{\al,s,q}$ }\label{sec:final_comp}
We are now concerned with the second of the problems listed above. That is, we solve the indeterminacy 
due to $\al=0$ in $\E_\al(0,q)$ and $\E_\al(1,q)$.

From the estimate
$$
1-(1-q)^\al>-\al\log(1-q)-\frac{\al^2}{2}\log^2(1-q),\qquad (\al,q)\in(0,1)\times(0,1)
$$
it follows that
$$
\E_\al(0,q)<\frac{2-q}{q}\left( -1+\frac{q(1-q)^{\al-1}}{-\log(1-q)-\frac{\al}{2}\log^2(1-q)}\right),
$$
$$
\E_\al(1,q)>\frac{2-q}{q}\left( 1-\frac{q}{-\log(1-q)-\frac{\al}{2}\log^2(1-q)}\right).
$$
Denote by $\tilde \E_\al(0,q)$ and $\tilde \E_\al(1,q)$ the rhs of the previous relations and define 
$$
Z_\al^f(s,q)\bydef \widetilde R_\al(s,q)(1+\tilde \E_\al(1,q))^2+\widetilde R_\al(1-s,q)(1+\tilde\E_\al(0,q))^2.
$$
Clearly, for those $(\al,s,q)\in\Lambda^1_{\al,s,q}$ where $\widetilde R_\al(1-s,q)$ is negative, we have $\widetilde{\mathcal I}^f_\al(s,q)>Z^f_\al(s,q)$. 
\\

To summarise:

 to conclude  that $\mathcal I^f_\al(s,q)>0$ for any $(a,q,s)\in \Lambda_{\al,s,q}$ we perform the following computations:
\begin{itemize}
\item $M^f_\al(s,q)>0$  in $\Lambda^0_{\al,s,q}$, refer to Appendix \ref{sec:AppMaf}
\item $Z^f_\al(s,q)>0$ in $\Lambda^1_{\al,s,q}$, refer to Appendix \ref{sec:AppZaf}
\item $ \widetilde {\mathcal I}_\al^f(s,q)>0$ in $\bar \Lambda_{\al,s,q}\setminus\{\Lambda^0\cup \Lambda^1\}=[\al^0,1]\times[\bar q,0.9]\times[0,\frac{1}{2}]$, refer to Appendix \ref{sec:AppIaftilde}
\end{itemize}

 \noindent {\it Proof of  Theorem \ref{theo:main}}

Let $\bar q=0.4 $ and $\alpha^0=0.4$. By rigorous computation it results that the quantities $M^f_\al(s,q)$, $Z^f_\al(s,q)$ and $\widetilde {\mathcal I}_\al^f(s,q)$ are all positive in $\Lambda^0_{\al,s,q}$, $\Lambda^1_{\al,s,q}$, $\bar \Lambda_{\al,s,q}\setminus\{\Lambda^0\cup \Lambda^1\}$ respectively. It follows that $\mathcal I^f_\al(s,q)>0$ for any $(a,q,s)\in \Lambda_{\al,s,q}$. The latter,  Proposition \ref{theo:infinite} and Proposition \ref{prop:0.9} provide the proof of Proposition \ref{prop:main}.
According to remark \ref{rmk:proof}, the proof of the Theorem \ref{theo:main} follows.
\qed

\section{Appendix}
All the rigorous computations have been performed in Matlab with the interval arithmetic package INTLAB \cite{Ru99a}.

For each of the rigorous computation reported in the paper, in the following we list the choice of the parameters and the subdivision of the domain. For a interval $A=[a_{min},a_{max}]$, to be decomposed into the union of subintervals $A\subset_i \cup A_i$, we refer to $a$-grid  as the set of points $\{a_j\}$ generating the subintervals $A_i=[a_i,a_{i+1}]$. When $\Delta a$ is a number it means that $a_i=a_{min}+i\Delta a$. In case $\Delta a$ is $equispaced(n)$, it means that $a_0=a_{min}$, $a_n=a_{max}$ and $a_j$ for $ 0\leq j \leq n$ are set numerically  equispaced.

\subsection{Checking $\partial_q \E(\al,q,1-s)>0$}
\subsubsection{Computation of $N_0(\al,q)>0$}\label{sec:AppN0} 
Let $\bar\al=0.15$ be fixed.
We check that $N_0(\al,q)$ given in \eqref{eq:N_0} is positive for any $\al\in[0,\bar\al]$ and $q\in[0.9,q_{max}]$ and that $\tilde N_0(\al,q)$ given in \eqref{eq:tildeN0} is positive for any $\al\in[0,\bar\al]$ and $q\in[q_{max},1]$.

For the first computation we  split the domain according to the following grid points 
	$$
	\begin{array}{cc}
	{\bf \al-grid} &\qquad  {\bf q-grid}\\
	\begin{array}{cl}
	& \Delta \al\\
	\hline
	0\leq \al\leq 0.04&0.002\\
	0.04\leq\al\leq \bar\al& 0.01
	\end{array}
	\quad 
	&
	\qquad \qquad \begin{array}{cc}
	& \Delta q\\
	\hline
	0.9\leq q \leq q_{max}&{\rm equispaced}(50)\\
	\\
	\end{array}
	\end{array}
	$$
and  results in $N_0(\al,q)\leq 0.0056$.
For the computation of $\tilde N_0(\al,q)$ we consider the grid produced by
$$
\begin{array}{cc}
	{\bf \al-grid} &\qquad  {\bf q-grid}\\
	\begin{array}{cl}
	& \Delta \al\\
	\hline
	0\leq\al\leq \bar\al& 0.01
	\end{array}
	\quad 
	&
	\qquad \qquad \begin{array}{cc}
	& \Delta q\\
	\hline
	q_{max}\leq q \leq 1&{\rm equispaced}(5)\\
	\end{array}
	\end{array}
$$
giving $\tilde N_0(\al,q)\geq 0.2721$.

\subsubsection{Computation of $N_1(\al,q)>0$}\label{sec:AppN1} 
Let $\hat\al=0.8$ be chosen.
We compute the function $N_1(\al,q)$ given in \eqref{eq:N1} for any $\al\in[\hat\al,1]$ and $q\in[0.9,\hat q]$ and the function $\tilde N_1(\al,q)$ \eqref{eq:N1tilde} for $(\al,q)\in[\hat \al,1]\times[\hat q,1]$. We check that both the functions are strictly positive in their own domain.
For the computation of $N_1(\al.q)$, the grid 
$$
	\begin{array}{cc}
	{\bf \al-grid} &\qquad  {\bf q-grid}\\
	\begin{array}{cl}
	& \Delta \al\\
	\hline
	\hat \al\leq \al\leq 0.9&0.002\\
	0.9\leq\al\leq 1& 0.01
	\end{array}
	\quad 
	&
	\qquad \qquad \begin{array}{cc}
	& \Delta q\\
	\hline
	0.9\leq q \leq \hat q&{\rm equispaced}(30)\\
	\\
	\end{array}
	\end{array}
	$$
	returns $N_1(\al,q)>0.01246$.
	
To check that $\tilde N_1(\al,q)$ is positive we use the grid
$$
\begin{array}{cc}
	{\bf \al-grid} &\qquad  {\bf q-grid}\\
	\begin{array}{cl}
	& \Delta \al\\
	\hline
	\hat \al\leq\al\leq 1& 0.04
	\end{array}
	\quad 
	&
	\qquad \qquad \begin{array}{cc}
	& \Delta q\\
	\hline
	\hat q\leq q \leq 1&{\rm equispaced}(5)\\
	\end{array}
	\end{array}
$$
and it results in  $\tilde N_1(\al,q)>1.7509$ in the domain of interest.

\subsubsection{Computation of $N(\al,q)>0$ for $\bar\al\leq \al\leq \hat\al$}\label{sec:AppN} 
Rather than computing $\tilde N(\al,q)$ given in \eqref{eq:Ntilde}, we provide the computation for 
$\frac{q}{\al^2}\tilde N(\al,q)$. It results in $\frac{q}{\al^2}\tilde N(\al,q)>0.0108$ for any $(\al,q)\in[\bar\al,\hat\al]\times[0.9,1]$. The computation has been performed  according to the grid points  
$$
\begin{array}{cc}
	{\bf \al-grid} &\qquad  {\bf q-grid}\\
	\begin{array}{cl}
	& \Delta \al\\
	\hline
	\bar \al\leq\al\leq \hat\al& 0.001\\
	\\
	\end{array}
	\quad 
	&
	\qquad \qquad \begin{array}{cc}
	& \Delta q\\
	\hline
	0.9 \leq q \leq 0.95& 0.0015\\
	0.95\leq q\leq 1 & 0.003
	\end{array}
	\end{array}
$$

\subsection{Computation of $\mathcal I_\al^f(s,q)$}
In this section we collect the details of the rigorous computation of $\mathcal I_\al^f(s,q)$ for $(\al,s,q)\in \bar\Lambda_{\al,s,q}$. As listed at the end of section \ref{sec:final_comp}, the proof of $\mathcal I_\al^f(s,q)>0$ amounts to providing that $M^f_\al(s,q)$, $Z^f_\al(s,q)$ and $\widetilde {\mathcal I}^f_\al(s,q)$ are all positive in their own domain of definition.
Let us fix $\bar q=0.4$.

\subsubsection{Computation of $M^f_\al(s,q)>0$ in $\Lambda^0_{\al,s,q}$}\label{sec:AppMaf}
 Note that  the function $M^f_\al(s,q)$ involves the quantities $K^\al_n(s)$ that, according to the definition \eqref{eq:Kn}, are not well defined at $s=0$. Therefore when computing $K^\al_n(S)$ with   any set $S$ that intersects the subspace $\{s=0\}$, we have to use the equivalent formulation \eqref{eq:K_nbis}. Moreover, it is advised to avoid  the division by intervals whose elements are close to zero, because that enlarges the size of the resulting interval.  Hence, we use the formula \eqref{eq:K_nbis} to compute the value $K^\al_n(s)$ for any $s\in(0,\frac{1}{2})$.
 
 The computation has been performed according to the grid
 $$
\begin{array}{ccc}
	{\bf \al-grid} &\qquad  {\bf q-grid}& \qquad {\bf s-grid}\\
	\begin{array}{cl}
	& \Delta \al\\
	\hline
	0\leq\al\leq 1& 0.025\\
	\end{array}
	\quad 
	&
	\qquad \qquad \begin{array}{cc}
	& \Delta q\\
	\hline
	0 \leq q \leq \bar q& 0.01\\
	\end{array}
	&
	\qquad \qquad \begin{array}{cc}
	& \Delta s\\
	\hline
	0 \leq s \leq 0.5& 0.01\\
	\end{array}
	\end{array}
$$
It results in $M_\al^f(s,q)>0.1223$.
\subsubsection{Computation of $Z^f_\al(s,q)>0$ in $\Lambda^1_{\al,s,q}$}\label{sec:AppZaf}
Choose $\al_0=0.4$. 
We first factor out the denominator $q^2(\log(1-q)+\frac{\al}{2}\log^2(1-q))^2$ and compute the remaining factor, still denoted by $Z^f_\al(s,q)$.
With the grid

 $$
\begin{array}{ccc}
	{\bf \al-grid} &\qquad  {\bf q-grid}& \qquad {\bf s-grid}\\
	\begin{array}{cl}
	& \Delta \al\\
	\hline
	0\leq\al\leq \al_0& 0.002\\
	\end{array}
	\quad 
	&
	\qquad \qquad \begin{array}{cc}
	& \Delta q\\
	\hline
	0 \leq q \leq \bar q& 0.001\\
	\end{array}
	&
	\qquad \qquad \begin{array}{cc}
	& \Delta s\\
	\hline
	0 \leq s \leq 0.5& 0.002\\
	\end{array}
	\end{array}
$$
The  result is $Z_\al^f(s,q)>0.01987$.

\subsubsection{Computation of $\mathcal{\widetilde I}_\al^f(s,q)>0$ in $\bar \Lambda_{\al,s,q}\setminus\{\Lambda^0\cup \Lambda^1\}$ }\label{sec:AppIaftilde}

The quantity \eqref{eq:tildeIfasq} is rigorously computed on the grid produced by the mesh
$$
\begin{array}{ccc}
	{\bf \al-grid} &\qquad  {\bf q-grid}& \qquad {\bf s-grid}\\
	\begin{array}{cl}
	& \Delta \al\\
	\hline
	\al_0\leq\al\leq 1& 0.005\\
	\end{array}
	\quad 
	&
	\qquad \qquad \begin{array}{cc}
	& \Delta q\\
	\hline
	\bar q \leq q \leq 0.9& 0.005\\
	\end{array}
	&
	\qquad \qquad \begin{array}{cc}
	& \Delta s\\
	\hline
	0 \leq s \leq 0.5& 0.005\\
	\end{array}
	\end{array}
$$
The response of the computation is $\widetilde {\mathcal I}^f_\al(s,q)>0.09260$ for any $(\al,s,q)\in\bar \Lambda_{\al,s,q}\setminus\{\Lambda^0\cup \Lambda^1\}$.

\bibliographystyle{model1-num-names}

\begin{thebibliography}{10}

\bibitem{LeVerrier1859}
U.~J. {Le Verrier}.
\newblock Theorie du mouvement de Mercure.
\newblock {\em Annales de l'Observatoire de Paris}, 5(1), 1859.

\bibitem{Newcomb}
S.~Newcomb.
\newblock Discussion and results of observations on transits of Mercury from
  1677 to 1881.
\newblock {\em Astronomical papers prepared for the use of the American
  ephemeris and nautical almanac}, 1:363--487, 1882.

\bibitem{Hall1894}
A.~Hall.
\newblock A suggestion in the theory of Mercury.
\newblock {\em Astronomical Journal}, 14(319):49--51, 1894.

\bibitem{Bertrand}
J.~Bertrand.
\newblock Th{\'e}or{\'e}me relatif au movement d'un point attir{\'e} vers un
  center fix{\'e}.
\newblock {\em Compt. Rend. Acad. Sci.}, 77:849--853, 1873.

\bibitem{Mcgehee}
R. McGehee.
\newblock Double collisions for a classical particle system with
  nongravitational interactions.
\newblock {\em Comment. Math. Helvetici}, 56:524--557, 1981.

\bibitem{gbegfugfg}
G. Bellettini, G. Fusco, and G.~F. Gronchi.
\newblock Regularization of the two-body problem via smoothing the potential.
\newblock {\em Comm. Pure Appl. Anal.}, 2:323--353, 2003.

\bibitem{Stoica}
C. Stoica and A. Font.
\newblock Global dynamics in the singular logarithmic potential.
\newblock {\em Journal of Physics A: Mathematical and General}, 36(28):7693,
  2003.

\bibitem{LeviCivita}
T.~Levi-Civita.
\newblock Sur la r{\'e}gularisation du probl{\`e}me des trois corps.
\newblock {\em Acta Mathematica}, 42(1):99--144, 1920.

\bibitem{MR2836348}
R. Castelli and S. Terracini.
\newblock On the regularization of the collision solutions of the one-center
  problem with weak forces.
\newblock {\em Discr. Contin. Dyn. Syst. Ser. A}, 31(4):1197--1218, 2011.

\bibitem{MR3057157}
R. Castelli, F. Paparella, and A. Portaluri.
\newblock Singular dynamics under a weak potential on a sphere.
\newblock {\em Nonlinear Differential Equations and Applications NoDEA},
  20(3):845--872, 2013.

\bibitem{MR2299759}
S. Terracini and A. Venturelli.
\newblock Symmetric trajectories for the 2N-body problem with equal masses.
\newblock {\em Archive for Rational Mechanics and Analysis}, 184(3):465--493,
  2007.

\bibitem{Richstone1982}
D.~O. Richstone.
\newblock Scale-free models of galaxies. II - A complete survey of orbits.
\newblock {\em Astrophysical Journal,}, 252:496--507, 1982.

\bibitem{Schwa1979}
M. Schwarzschild.
\newblock A numerical model for a triaxial stellar system in dynamical
  equilibrium.
\newblock {\em Astrophys. Journal}, 232:236--247, 1979.

\bibitem{Newton}
P. Newton.
\newblock {\em The $N$-vortex problem. Analytical techniques.}
\newblock Number 145 in Applied Mathematical Sciences. Springer-Verlag, New
  York, 2001.

\bibitem{MR2029363}
C.~E. Kenig, G. Ponce, and L. Vega.
\newblock On the interaction of nearly parallel vortex filaments.
\newblock {\em Communications in Mathematical Physics}, 243(3):471--483, 2003.

\bibitem{MNR:MNR22071}
S.~R. Valluri, P.~A. Wiegert, J.~Drozd, and M.~Da~Silva.
\newblock A study of the orbits of the logarithmic potential for galaxies.
\newblock {\em Mon. Not. R. Astron. Soc.}, 427(3):2392--2400, 2012.

\bibitem{Valsecchi}
F.~{Valsecchi}, W.~M. {Farr}, B.~{Willems}, C.~J. {Deloye}, and V.~{Kalogera}.
\newblock Tidally induced apsidal precession in double white dwarfs: A new mass
  measurement tool with LISA.
\newblock {\em The Astrophysical Journal}, 745(2):137, 2012.

\bibitem{Feldman1979}
G. Feldman, T. Fulton, and A. Devoto.
\newblock Energy levels and level ordering in the wkb approximation.
\newblock {\em Nuclear Physics B}, 154(3):441--462, 8 1979.

\bibitem{valluri99}
S.~R. Valluri, R.~Biggs, W.~Harper, and C.~Wilson.
\newblock The significance of the Mathieu-Hill differential equation for
  Newton's apsidal precession theorem.
\newblock {\em Canadian Journal of Physics}, 77(5):393--407, 1999.

\bibitem{Danby}
J.~M. Danby.
\newblock {\em Fundamentals of Celestial Mechanics}.
\newblock Willmann-Bell, Inc., 2nd ed., rev. ed. Richmond edition, 1988.

\bibitem{whitta}
E.~T. Whittaker.
\newblock {\em A Treatise on the Analytical Dynamics of Particles and Rigid
  Bodies, 4th Ed.}
\newblock Cambridge University Press, Cambridge, 1959.

\bibitem{MNR:MNR8819}
S.~R. Valluri, P.~Yu, G.~E. Smith, and P.~A. Wiegert.
\newblock An extension of newton's apsidal precession theorem.
\newblock {\em Mon. Not. R. Astron. Soc.}, 358(4):1273--1284, 2005.

\bibitem{1997MNRAS}
J.~Touma and S.~Tremaine.
\newblock A map for eccentric orbits in non-axisymmetric potentials.
\newblock {\em Mon. Not. R. Astron. Soc.}, 292:905--919, 1997.

\bibitem{MR2497199}
F.~C. Santos, V. Soares, and Al.~C. Tort.
\newblock Determination of the apsidal angles and Bertrand's theorem.
\newblock {\em Phys. Rev. E}, 79(3):036605, 2009.

\bibitem{struck}
C.~Struck.
\newblock Simple, accurate, approximate orbits in the logarithmic and a range
  of power-law galactic potentials.
\newblock {\em The Astronomical Journal}, 131(3):1347, 2006.

\bibitem{MR1433938}
S.~R. Valluri, C. Wilson, and W. Harper.
\newblock Newton's apsidal precession theorem and eccentric orbits.
\newblock {\em Journal for the History of Astronomy}, 13, 1997.

\bibitem{MR3159801}
R.~Castelli.
\newblock A study of the apsidal angle and a proof of monotonicity in the
  logarithmic potential case.
\newblock {\em Journal of Mathematical Analysis and Applications}, 413(2):727
  -- 751, 2014.

\bibitem{LambertW}
R.~M. Corless, G.~H. Gonnet, D.~E.~G. Hare, D.~J. Jeffrey, and D.~E. Knuth.
\newblock On the Lambert W function.
\newblock {\em Adv. Comput. Math.}, 5:329--359, 1996.


\bibitem{grant}
A.~K. Grant and J.~L. Rosner.
\newblock Classical orbits in power-law potentials.
\newblock {\em Am. J. Phys}, 62(4):310--315, 1994.

\bibitem{MR2652784}
S.~M. Rump.
\newblock Verification methods: rigorous results using floating-point
  arithmetic.
\newblock {\em Acta Numer.}, 19:287--449, 2010.

\bibitem{intervalbook}
R.~E. Moore, R.~B. Kearfott, and M.~J. Cloud.
\newblock {\em Introduction to Interval Analysis}.
\newblock Society for Industrial and Applied Mathematics, 2009.

\bibitem{Ru99a}
{S.~M.} Rump.
\newblock {INTLAB - INTerval LABoratory}.
\newblock In Tibor Csendes, editor, {\em {Developments~in~Reliable Computing}},
  pages 77--104. Kluwer Academic Publishers, Dordrecht, 1999.
\newblock http://www.ti3.tu-harburg.de/rump/.

\end{thebibliography}

\end{document}